\documentclass[a4paper,12pt]{amsart}\usepackage{amsfonts,amssymb,bbm,enumerate,mathrsfs,abstract}
\usepackage{color}
\usepackage[bookmarksdepth=3]{hyperref}

\UseRawInputEncoding

 \def\C{\mathbb{C}}
\def\k{\Bbbk}
\def\K{\mathbbm{K}}\def\N{\mathbb{N}}\def\Q{\mathbb{Q}}
\def\R{\mathbb{R}}\def\Z{\mathbb{Z}}
\def\bk{{\bar{\k}}}

\def\di{\partial}

\def\bl{\langle}\def\br{\rangle}

\def\liml{\lim\limits}

\DeclareMathOperator{\Der}{Der}\DeclareMathOperator{\Frac}{Frac}
\DeclareMathOperator{\ord}{ord}

\DeclareMathOperator{\Span}{Span}\DeclareMathOperator{\Spec}{Spec}
\DeclareMathOperator{\Supp}{Supp}

\newcommand{\quot}[2]{{\footnotesize\left.\raisebox{1.2ex}{$#1$}\!\! \ensuremath\diagup \!\!\raisebox{-1.2ex}{$#2$}\right.}}
\newcommand{\quots}[2]{{\footnotesize\left.\raisebox{0.4ex}{$#1$}\! / \!\raisebox{-0.4ex}{$#2$}\right.}}

 \renewcommand{\stackrel}[2]{\ \lower 0.4ex \hbox{$\mathrel{\mathop{#2}\limits^{\scriptscriptstyle {#1}}}$}\ }

\def\ta{\tilde{a}}\def\tc{\tilde{c}}\def\td{{\tilde{d}}}\def\tf{{\tilde{f}}}

\def\tx{{\tilde{x}}}\def\txi{{\tilde{\xi}}}
\def\ty{{\tilde{y}}}\def\tz{{\tilde{z}}}

\def\hx{\hat{x}}

\def\al{\alpha}\def\de{\delta}\def\De{\Delta}
\def\ep{\epsilon}
\def\la{\lambda}\def\om{\omega}\def\si{\sigma}

\def\cm{{\frak m}}
\def\cU{\mathcal U}

\def\um{{\underline{m}}}
\def\ux{\underline{x}}\def\uy{{\underline{y}}}

\newcommand{\ber}{\begin{array}{l}}\newcommand{\eer}{\end{array}}
\newcommand{\bpm}{\begin{pmatrix}}\newcommand{\epm}{\end{pmatrix}}
\newcommand{\bbm}{\begin{bmatrix}}\newcommand{\ebm}{\end{bmatrix}}
\newcommand{\bM}{\begin{matrix}}\newcommand{\eM}{\end{matrix}}
\newcommand{\bee}{\begin{enumerate}}\newcommand{\eee}{\end{enumerate}}
\newcommand{\bei}{\begin{itemize}}\newcommand{\eei}{\end{itemize}}

\def\sset{\!\subset\!}\def\sseteq{\!\subseteq\!}\def\ssetneq{\!\subsetneq\!}\def\smin{\!\setminus\!}
\def\iff{if and only if }

\newtheorem{Statement}{Statement}[section]

\newtheorem{Theorem}[Statement]{Theorem}\newcommand{\bthe}{\begin{Theorem}}\newcommand{\ethe}{\end{Theorem}}
\newtheorem{Lemma}[Statement]{Lemma}\newcommand{\bel}{\begin{Lemma}}\newcommand{\eel}{\end{Lemma}}
\newtheorem{Proposition}[Statement]{Proposition}\newcommand{\bprop}{\begin{Proposition}}\newcommand{\eprop}{\end{Proposition}}
\newtheorem{Corollary}[Statement]{Corollary}\newcommand{\bcor}{\begin{Corollary}}\newcommand{\ecor}{\end{Corollary}}
\newtheorem{Definition}[Statement]{Definition}\newcommand{\bed}{\begin{Definition}}\newcommand{\eed}{\end{Definition}}
\newtheorem{Definition-Proposition}[Statement]{Definition-Proposition}

\def\bpr{~\\{\em Proof.\ }}
\newcommand{\epr}{{\hfill\ensuremath\blacksquare}}
\newtheorem{Remark}[Statement]{Remark}\newcommand{\beR}{\begin{Remark}\rm}\newcommand{\eeR}{\end{Remark}}
\newtheorem{Example}[Statement]{Example}\newcommand{\bex}{\begin{Example}\rm}\newcommand{\eex}{\end{Example}}

\newcommand{\bet}{\begin{tabular}{cccccccc}}\newcommand{\eet}{\end{tabular}}
\newcommand{\beq}{\begin{equation}}\newcommand{\eeq}{\end{equation}}

\def\into{{\hookrightarrow}}

\newcommand{\bin}[2]{\binom{#1}{#2}}
 \newcommand{\isom}[1]{\xrightarrow[\,\smash{\raisebox{1.15ex}{\ensuremath{\scriptstyle\sim}}}\,]{#1}}

\title[]{W\MakeLowercase{hen does a derivation of a ring admit the exponential?}}
\author[]{G\MakeLowercase{enrich} B\MakeLowercase{elitskii,}
A\MakeLowercase{lberto} F. B\MakeLowercase{oix and}
D\MakeLowercase{mitry} K\MakeLowercase{erner}}
\address{Department of Mathematics, Ben Gurion University of the Negev, P.O.B. 653, Be'er Sheva 84105, Israel.
 Department of Mathematics, Universitat Polit\`ecnica de Catalunya BarcelonaTech, Av. Eduard Maristany 16, 08019, Barcelona, Spain.}
\email{\!\!genrich@math.bgu.ac.il, \!alberto.fernandez.boix@upc.edu, \!dmitry.kerner@gmail.com}

\date{\today.\ \  Filename: \jobname.tex}

\thanks{A. F. Boix was partially supported by  Spanish Ministerio de Ciencia e Innovaci\'on grant PID2022-137283NB-C22}
\thanks{D. Kerner  was supported by the Israel Science Foundation,  grants No.  1910/18 and 1405/22}

\begin{document}

  \maketitle

 \begin{abstract}  
 Exponentials of (real/complex) vector fields are classically defined via  the vector field integration.
   Take a $\k$-algebra $\k[x] \sset R\sset \k[\![x]\!],$ where  $\k\supseteq \Q$ is a local domain.
   Suppose a derivation $\xi$ is $x$-adically nilpotent.
     Define the exp-operator via the Taylor expansion, $e^\xi:=\sum \frac{\xi^j}{j!}.$ It is a formal automorphism,
      $e^\xi\in Aut_\k(\k[\![x]\!]).$
   When does $e^\xi$ act on $R?$ When does the formal power series $e^\xi x\in \k[\![x]\!]$ belong to $R?$

\medskip

   We address this question for the following rings.
   \bee[i.]
   \item The algebraic power series, $R=\k\bl x\br,$   differentially finite (holonomic) power series, $D(\k[x]),$
    and their higher versions,   Picard-Vessiot extensions   $D^\bullet(\k[x]),$  Picard-Vessiot closure $D^\infty(\k[x]),$
      and  differentially-algebraic power series $D^{alg}(\k[x]).$
   \item Power series over normed fields. In particular, power series with coefficients of controlled growth, e.g. analytic/Denjoy-Carleman/Gevrey classes.
   \item Germs of smooth functions $\quots{C^\infty(\R^n,o)}{J},$ for arbitrary ideal $J\sset C^\infty(\R^n,o).$
   \eee

\medskip

In case i.   the operator $e^\xi$ is transcendental, and the power series $e^\xi x$ is ``usually" far from being algebraic.
 We give various criteria on  $e^\xi x$ to belong to $\k\bl x\br,$ $D(\k[x]),$  $D(R),$ or $D^{alg}(\k[x]).$
  Among our results are:
 \bei
 \item  the ring $D^\infty(\k[x,y])$ does not admit the Weierstrass preparation theorem.
\item (one variable case) If $e^\xi x\in D^\infty(R),$  then $e^\xi x\in R^{hens}$ (the Henselization)
\item (one variable polynomial case)  $e^\xi x\in D^\infty(\k[x])$ \iff $\xi =c\cdot x^p\cdot \di_x,$ for some $c\in \k.$
\eei

\medskip

In case ii. the answer is positive (i.e. $e^\xi$ acts on $R$) under rather weak assumptions on $R.$

In case iii. the answer is ``totally negative". For any $\xi\neq0 $ the operator $e^\xi$ (defined as before) does not act on the quotients of the ring of germs of smooth functions, $\quots{C^\infty(\R^n,o)}{J}.$
  \end{abstract}

\setcounter{secnumdepth}{6} \setcounter{tocdepth}{1}\tableofcontents

\section{Introduction}
\subsection{}Let $\k\!\supseteq\! \Q$ be a local domain, e.g. a field. Let $R$ be a
  (not necessarily Noetherian) $\k$-algebra satisfying: $\k[x]\sset R\sseteq \k[\![x]\!].$
 We use multi-variables, $x=(x_1,\dots,x_n)$.

Take a $\k$-linear derivation  $\xi\in \Der_\k(R).$
For various purposes one needs ``to integrate $\xi$ to a time-one flow", i.e.
to assign to $\xi$ the ($\k$-linear) automorphism of the algebra, $e^\xi \in
Aut_\k(R)$.
 In Differential Geometry (for $\k=\R,\C$) the operator $e^\xi $ is defined by  vector field integration, resolving the system of ODE's, $x'=\xi(x).$
 For   $\k\supseteq\Q$ a local domain (with no norm)  the natural definition is via the Taylor expansion,
$e^\xi\!:=\!\sum^\infty_{j=0}\frac{\xi^j}{j!}.$   The  old  question reads: {\em  When  does  $e^\xi$  act  on  $R$?}

\medskip

If $R=\k[x]$ and $\xi$ is locally nilpotent\footnote{i.e. for each $f\in R$ there exists $N\in \N$ ensuring $\xi^N f=0$},
 then one gets a well defined automorphism $e^\xi\circlearrowright R.$ These locally nilpotent derivations have
  been studied extensively, see e.g. \cite{Freudenburg}.

\medskip

In local Geometry/Algebra, Singularity Theory, Dynamical Systems (and neighboring areas) one works with   local $\k$-algebras,
 $\k[x]_{(x)}\sseteq R\sseteq \k[\![x]\!].$  Then {\em no  (non-zero) derivation is locally nilpotent,} see \S\ref{Sec.Prelim.Derivations}.
If $\k$   has no pre-given topology, then even the coefficients of monomials in the power series $e^\xi  f$ are not well defined. (They contain infinite summations.)
  Hence one restricts to ``$x$-adically nilpotent" derivations,
 \beq\label{Eq.x-adically.nilpotent.derivations}
 \Der^{nilp}_\k(R):=\{\xi \ | \  \xi^N(x)\sseteq (x)^2  \ for\ N\gg1\}\sset  \Der_\k(R).
 \eeq
 (And then $\xi^N(x)^d\sseteq (x)^{d+1}$ for each $d\in \N$ and a corresponding $N\gg1$.)
 Geometrically these are vector fields on $(\k^n,o)$ that vanish at the origin and whose linear part is nilpotent.

\medskip

  This $x$-adic nilpotence ensures (by the $x$-adic completeness): $e^\xi$ is a well-defined (formal) automorphism of the ring $\k[\![x]\!],$
    \cite[Chapter II, \S6]{Bourbaki.Lie}.
   Thus the initial question becomes:
  {\em Does the action $e^\xi\circlearrowright \k[\![x]\!]$ restrict to the action $e^\xi\circlearrowright R?$}
 The classical  answers are:
\bei
\item    ``yes" for the analytic ring $\k\{x\},$  archimedean or non-archimedean,      \cite[Chapter II,\S7,\S8]{Bourbaki.Lie};
 ``yes" for Gevrey classes, \cite{BrocheroLopezHernanz2009};
\item  ``almost never" for the ring of algebraic power series, $\k\bl x\br.$  The power series $e^\xi x$ is typically non-algebraic, and was
 considered ``completely transcendental".
\eei
For other rings, $\k[x]_{(x)}\sset R\sset\k[\![x]\!],$ only some particular (positive and negative) examples are  known.

\medskip

If $\xi$ acts on $R,$ then, in particular, $e^\xi x\in R.$
 The converse holds for rings admitting compositions, because $e^\xi f(x)=f(e^\xi x)$.
   Thus the essential question is whether $e^\xi x\in R,$ or ``How far is the Lie series $e^\xi x$ from $R$?"

For the rings like $\k\bl x\br,$ $D^\bullet(\k[x])$ this goes in the style of Liouville theory/Differential Galois theory,
see e.g. \cite{Raab-Singer} and \cite{van der Put-Singer}.
 But the involved ODE, $x'=\xi(x),$ is non-linear. Thus the standard known results do not seem to be directly applicable.

\subsection{}\label{Sec.1.2} In the formal/analytic case the action $e^\xi\circlearrowright R$ gives the Lie map
  $\Der^{nilp}_\k(R)\stackrel{exp}{\to}Aut_\k(R).$
 This allows to study the neighborhood of the unit element, $Id\in Aut_\k(R)$, via the ``tangent space to the group", $\Der^{nilp}_\k(R)$,
  \cite[Chapter III]{Bourbaki.Lie}. Hence the importance of the exponentiability.
  More generally, these ``Lie series of vector fields" are widespread in Mathematics and Physics, see e.g. \cite{Olver}, \cite{Milnor}, \cite{Winkel}. In the analytic case they have been studied exhaustively through 20'th century.

\medskip

In many cases $e^\xi x\not\in R,$ thus  the map $exp$ is not  defined on the whole module $\Der^{nilp}_\k(R).$
 Then one asks more specific questions, dictated by applications.
 \bei
 \item
On which (largest) subset of $\Der^{nilp}_\k(R)$ is the map $exp$ defined?
 We call $\xi$ {\em exponentiable over $R$} if $e^\xi x\in R.$
 Denote the subset of exponentiable derivations by $\Der^{exp}_\k(R)\sseteq \Der^{nilp}_\k(R).$

 A remark, if $e^\xi$ acts on $R$ then so does $e^{c\cdot \xi}$ for any $c\in \k,$ see Lemma \ref{Thm.TFAE}.
  Thus the subset $\Der^{exp}_\k(R)\sseteq \Der^{nilp}_\k(R)$  is a $\k$-cone. But it is often not a $\k$-submodule.

 \item (The complexity of $e^\xi x$ for a given $\xi$) Is $e^\xi x$ contained in a prescribed extension of $R?$
 E.g. for $D$-finite power series (and their higher analogs), can one ensure $e^\xi(D^j(\k[ x]))\sseteq D^{j+d}(\k[ x])$ for some   $d$ and   $j?$
 Or (at least) $e^\xi(D^j(\k[ x]))\sseteq D^\infty(\k[ x])?$

     \item For many rings, even if $e^\xi$ does not act on $R,$ this operator is $(x)$-adically approximated by automorphisms of $R.$
 This is frequently used for the study of the group $Aut_\k(R)$ in terms of the module $\Der_\k(R).$
 A more delicate question is:
     which derivations $\xi\in   \Der^{nilp}_\k (R)$ can be $(x)$-adically approximated by exponentiable derivations?
 Namely: describe the $(x)$-adic closure  $\overline{\Der^{exp}_\k(R)}\sseteq \Der^{nilp}_\k(R).$

         Suppose $\xi$ is (non)exponentiable, are its perturbations by terms of high enough orders (non) exponentiable?
           Is the property of being (non)exponentiable finitely determined?

       \item How to compute $e^\xi x$ effectively? Recall that $e^\xi x$ is the ``time-one flow",
        $e^{t\xi}x|_{t=1},$ for the (system of) ODE $x'=\xi(x).$
        While the flow is transcendental in general, $e^\xi x$ has ``neat formula" in some particular cases.
       E.g.
\beq\label{Eq.exp.der.for.x^p}
\text{for} \      \xi=c\cdot x^{p+1} \frac{d}{d x } , \text{ with } c\in \k \text{ and } p\ge 1,
 \text{ one has: }
  e^\xi x= \frac{x}{\big(1-c\cdot p \cdot x^p\big)^\frac{1}{p}}\in \k\bl x\br.
  \eeq
 This classical formula is widely used in Mathematis and Physics, e.g. \cite{Milnor}, \cite{Olver}, \cite{Schottenloher}.
           More generally, one needs a closed formula for $e^\xi x$ in terms of some ``finite elementary operations",
           e.g. algebraic extensions,   antiderivatives, and resolving implicit function equations.
 \eei

We address these questions for the three (most important) classes of rings of the abstract.

\subsection{The results}
The standard starting point is: exponentiability vs solvability of ODE's.
\\\mbox{{\bf Lemma  \!\ref{Thm.TFAE} \!(roughly)}
   {\em $e^\xi x\!\in\! R$ \iff the ODE-system $x'\!=\!\xi(x)$ is solvable over $R.$}}

   Hence the initial question becomes: ``Is $R$ ODE-closed?"
    Or, for which systems of ODE's over $R$ the (unique, formal) solution belongs to $R?$

  This question was studied for some subrings of the ring $C^\infty(\cU)$ of smooth functions, see e.g. \cite{Rainer-Schindl}.
 Many results are known for Gevrey classes, see \cite{Cano} and references therein.

\subsubsection{The one-variable case,  $n=1,$ and the rings $\k\bl x\br,$ $D^j(\k[x]),$ $D^\infty(\k[x]),$ $D^{alg}(\k[x])$}\label{Sec.1.3.1}
 (See \S\ref{Sec.Prelim.Different.Extension} for the definitions.)
 Recall that $e^\xi x$ is the time-one flow of the vector field $\xi.$ Therefore $e^\xi x$ can be computed by resolving the ODE $x'=\xi(x).$
 This ODE is integrated (in the completely standard way) to get the implicit function equation:

{\bf Lemma \ref{Thm.Exp.Der.Case.n=1}.} {\em
1. The power series $y(x)=e^\xi x\in \k[\![x]\!]$  is determined by the    equation
\beq
  \frac{1}{c\cdot a(y)^{p}} - \frac{1}{c\cdot a(x)^{p}}  +res_o\om\cdot ln\frac{y}{x}=1.
  \eeq
(Here $res_o\om$ is the residue of the associated 1-form.)

 2.  In particular, if $res_o\om\!=\!0$  then  $e^\xi x\!=\!a^{-1}[\frac{a(x)}{\big(1+c\cdot a(x)^{p}   \big)^\frac{1}{p}}].$
}

 Here $a^{-1}[\dots]$ is the inverse of the power series $a(x)$.

Thus  $e^\xi x$ is not ``arbitrarily complicated",
 but belongs to a controlled extension of the initial ring $R.$ This  has numerous consequences.
\bei
\item
In several ``numerical cases" we have neat formulas, with $e^\xi x\in D^2(\k[x]),$ see  Example \ref{Ex.n=1.Computations}.
  But immediately we hit the nasty case:
\\
{\bf Example \ref{Ex.Lambert.Function}.} {\em $e^\xi x\not\in D^\infty(\k[x])$ for $\xi=-\frac{x^2 e^{-x}}{x-1}\di_x.$}
\\
This shows how complicated the power series $e^\xi x$ can be, even if we start from $\frac{x^2 e^{-x}}{1+x}\!\in\! D^1(\k[x]).$
 (Recall that $D^\infty(\k[x])$ is the Picard-Vessiot closure of $\k[x]$ inside $\k[\![x]\!]$.)

\item
An offshoot   of this nasty case is the pathology of the Picard-Vessiot closure:
\\
{\bf Corollary \ref{Thm.no.WPT.for.Dinfty}.}
 {\em The ring $D^\infty(\k[x,y]),$ being regular, local,  Henselian, does not admit the Weierstrass preparation theorem.}
\hfill (In particular, $D^\infty(\k[\ux])$ is not a W-system.)

This (unpleasant) result seems to be new.
 Examples of regular local Henselian rings with no Weierstrass preparation are known in Analysis e.g. Gevrey or
  quasi-analytic Denjoy-Carleman rings  for $n\ge2,$ \cite{A.B.B.N.Z.},  \cite{Parus.Rol.}.
 But the ring  $D^\infty(\k[x,y])$ is the first ``purely algebraic" case,  over arbitrary field $\k.$

\item
For the much larger ring of differentially algebraic power series the statement is positive (and probably known):
\\{\bf Corollary \ref{Thm.n=1.exp.for.Diff.algebraic}.} \ {\em 1. $e^\xi x\in D^{alg}(\k[x]) $ for any   $\xi\in \Der^{nilp}(D^{alg}(\k[x])).$
\\2. For algebraic derivations, $\xi\in  \Der_\k^{nilp}(\k\bl x\br),$ the differential order of $e^\xi x$ is 1.}

\item
Another application of Lemma \ref{Thm.Exp.Der.Case.n=1} is the density result:
\\{\bf Corollary \ref{Thm.exp.for.exact.form}.}
 {\em The $(x)$-adic closure of the set of exponentiable derivations, $\overline{\Der^{exp}_\k(\k\bl x\br)}\sset  \Der_\k^{nilp}(\k[\![ x]\!]),$
is defined  by the single condition, $res_o\om=0.$}

\eei

\medskip

\noindent The natural goal (in view of \S\ref{Sec.1.2}) is: assuming $e^\xi x\!\in \!D(R),$ for a subring $R\sset \k[\![x]\!],$
    to push $e^\xi x$ ``closer to $R$".
For   algebraic power series, $\xi\!\in  \! \Der^{nilp}_\k(\k\bl x\br),$ we get a  strong criterion:
\\
{\bf Theorem \ref{Thm.D^infty.to.Rhens}.} {\em If $e^\xi x \in D^\infty(R) $ then $e^\xi x \in R^{hens}.$}  (The Henselization of $R.$)
\\{\bf Corollary \ref{Thm.If.D^infty.then.algebraic}.}
 {\em In particular,   $e^\xi x\in D^\infty(\k[x])$ \iff $e^\xi x\in \k\bl x\br.$}

\medskip

This (completely new) Liouville-type dichotomy can be stated as ``the operator $e^\xi$ is orthogonal to the hierarchy  $D^1(R)\sset D^2(R)\sset\cdots$".

\medskip

 For polynomial derivations the result is even stronger, showing that \eqref{Eq.exp.der.for.x^p}
 is the only case with a  $D^\infty$-answer:
\\
{\bf Theorem \ref{Thm.exp.Der.is.D1.iff.xi=x^p}.} {\em
Let $\xi=a(x)\cdot \di_x$ with $a(x)\in (x)^2\sset \k[x].$
Then $e^\xi x\in D^\infty(\k[x])$ \iff $\xi =c\cdot x^p\cdot \di_x,$ for some $c\in \k.$
}

\subsubsection{The case of several variables, $n\ge2$}
Here  the study of $e^\xi x$ is as complicated as the solution of systems of (arbitrary, nonlinear) ODE's.
 Here enter the difficulties of first integrals.
  In \S\ref{Sec.Exp.Der.computation.Separated.Derivations}
 we give explicit
 formulas only for derivations with ``separable variables".
 (In this case one has immediate full system of first integrals.)

\subsubsection{Power  series in $n$ variables, with  coefficients  in a normed field.}  Here
    derivations are exponentiable  under rather   weak assumptions.
 Let $\k[x]\sset R\sset\k[\![x]\!],$ where $\k$ is a normed field.

 {\bf Theorem \ref{Thm.Exp.Der.Bounded.Coefficients}.} {\em Suppose $R$  is ``closed under majorization of coefficients",
  and satisfies the ``implicit function theorem in two variables".
  Then    $e^\xi x\in R$   for any $\xi\in   \Der^{nilp}_\k(R).$}

\medskip

\mbox{The proof goes by reduction to the one-variable case ($n\!=\!1$). Then we invoke Lemma \ref{Thm.Exp.Der.Case.n=1}.}

In the simplest case of analytic power series, $\k\{x\},$ we get the classical Cauchy theorem for analytic flows of
  vector fields (over arbitrary normed field, e.g. p-adics).
 Unlike  the numerous  proofs of that Cauchy theorem, e.g. see the references in \cite{Carrillo}, ours does not involve any bounds or estimates.

More interesting cases for Theorem \ref{Thm.Exp.Der.Bounded.Coefficients} are rings with coefficients of controlled growth,
 e.g. Denjoy-Carleman and Gevrey classes, Example \ref{Ex.Exp.Der.acts.for.normed.field}.

\subsubsection{Smooth function germs.} For the ring $R=\quots{C^\infty(\R^n,o)}{J},$ assuming the zero locus of $J$ is not just the origin $o\in \R^n,$
  the answer is totally negative.
\\{\bf Theorem \ref{Thm.C.infty.case.exp.does.not.act}.}
{\em  For any derivation $0\!\neq \!\xi\!\in   \!\Der_\R (R)$ the operator $e^\xi $  does not act on $R.$
}

Namely, for any given $\xi$ we construct a function-germ $f\in C^\infty(\R^n,o)$ such that $e^\xi f$
  diverges on a convergent sequence of points $o\leftarrow \{z_\bullet\}\sset V(J)$.

This shows (again) how complicated is the ring  $C^\infty(\R^n,o),$ as compared to various subrings of $\k[\![x]\!].$

 The proof is of (quite non-algebraic) $C^\infty$-nature.
 The statement is immediate when $\xi|_o\neq0.$ (Then $\xi$ can be rectified, and one invokes Borel's lemma.)
 But the case of $\xi|_o=0$ is completely non-trivial.

\subsubsection{The further work}
  These results are just the peak of the iceberg. The operator $e^\xi$ is the common/widespread tool
 in Mathematics and Physics, hence every property is  of immediate value.
\\
Our methods are applicable to various other rings, e.g. the traditional rings of Transcendence Theory
 and Diophantine Approximations (or Analytic Combinatorics),
  Mahler functions, G-functions, E-functions.
 Some of our results extend naturally to the case of  $n\ge2.$

 All this stuff is postponed to the subsequent paper.

\subsection{Acknowledgements}
We thank D. Novikov and G. Rond for useful advice.

 Most of the work was done during the stay of A. F. Boix in BGU, but finalizing the results was delayed by the covid-times and by the war.

\subsection{Declarations}
The author declare that they have no conflict of interest. There are no data sets generated during the current
study.

\subsection{Notations and conventions}\label{Sec.Prelim.Notations.Conventions}
\bee[\!\!\!\bf i.]
\item
We use   multi-variables,  $ x=(x_1,\dots,x_n),$    $ x^o=(x^o_1,\dots,x^o_n),$ and $x^m=x^{m_1}_1\cdots x^{m_n}_n.$

\item
Through the paper $\k\supseteq \Q$ is a local domain, the main example being a field.
 Occasionally we pass to the  field of fractions  $\Frac(\k).$  If $\k$ is a field, then $\Frac(\k)=\k.$

\item
  Except for \S\ref{Sec.Cinfty.case}, $R$ is  a differential (not necessarily Noetherian) $\k$-subalgebra,  $\k[x] \sset R\sseteq \k[\![x]\!] .$
 Thus the partials $\di_{x_1},\dots\di_{x_n}$ act on $R.$

    Denote by   $\Frac(R)$  the field of fractions.

\item
  One says ``$R$   admits  compositions"  if  $f(g_1(x),\dots,g_n(x))\!\in\! R $ for all $f(x)\!\in\! R$ and for all $g_1(x),\dots,g_n(x)\!\in\! (x)\sset  R.$

\eee

\section{Preliminaries}\label{Preliminaries}

\subsection{Derivations}\label{Sec.Prelim.Derivations}
\bee[\!\!\bf i.]

\item  Let $R=\k[\![x]\!].$
The module of $\k$-linear derivations is a free $R$-module,   $\Der_\k(R)=R\bl \di_{x_1},\dots,\di_{x_n} \br,$
  see Theorem 30.6(i) of \cite{Matsumura}.
\\
The $(x)$-adically nilpotent derivations, see \eqref{Eq.x-adically.nilpotent.derivations}, satisfy:
   $(x)^2\!\cdot\! \Der_\k(R) \sseteq  \Der_k^{nilp}(R) \sset   (x)\! \cdot \!\Der_\k(R).$

For $n=1$ one has the equality: $ \Der_k^{nilp}(R)= (x)^2\cdot \Der_\k(R).$

For $n\ge 2$ the subset $ \Der_k^{nilp}(R)\sset \Der_k (R)$ is not a $\k$-submodule.

For $\xi\in \Der_k^{nilp}(R)$ the power series $e^{t\xi}x\in \k[\![t,x]\!]$ belongs  to $\k[t][\![x]\!].$

 \item Fix a $\k$-subalgebra $\k[x]\sseteq R\sseteq \k[\![x]\!].$
  A derivation $\xi\in \Der_\k R$ is called locally nilpotent on $R$ if $\xi^N f=0$ for each $f\in R$ and a corresponding $N=N_f.$
\\
{\bf Claim:} {\em Any ring $\k[x]_{(x)}\sseteq R\sseteq \k[\![x]\!]$ has no  non-zero  locally nilpotent derivations.}
 \bpr Take a derivation $0\neq \xi\in \Der_\k(R),$ suppose it is locally nilpotent.  Take a polynomial $f\in (x)^2\sset \k[x] $ satisfying: $\xi^{p+1} f=0$ and $\xi^p f\neq0$ for some $p\ge1.$ Then $\xi^{d\cdot p}(f^d)\stackrel{up\ to\ \k^\times}{\sim} (\xi^{  p}f)^d\neq0,$ while $\xi^{d\cdot p+1}(f^d)=0.$
\\Take $ \!\sum^\infty_{j=0} f^j\!=\!\frac{1}{1-f}\!\in\! \k[x]_{(x)}.$ Then $\xi^{d\cdot p}(\sum^{d-1}_0 f^j)\!=\!0.$
   But $\xi^{d\cdot p}(  f^d)\!\neq\!0$ and $ord (\xi^{d\cdot p}(  f^d))\!<\!ord (\xi^{d\cdot p}(  f^{d+1})).$
    Therefore $\xi^{d\cdot p}( \!\sum^\infty_{j=0} f^j)\!\neq\!0$ for each $d\!\in\! \N.$ In contradiction to $\xi$ being locally nilpotent.
 \epr

\item Occasionally one   changes the coordinates on $(\k^n,o),$  i.e. takes ($\k$-linear) automorphisms $\Phi\in Aut_\k R,$ $x\to \phi(x)=\tx(x).$
 This coordinate change acts on the derivation $\xi=\sum \xi_i \di_{x_i}$ by pushforward, $\phi_* \xi=\sum \xi_i \frac{\di \tx_j}{\di x_i }\di_{\tx_j}.$ We get: $\phi (e^\xi f(x))=e^{\phi_* \xi}f(x(\tx)).$
  In particular, the operator $e^\xi$ acts on $R\sset \k[\![x]\!]$ \iff the operator $e^{\phi_* \xi}$ acts on $\phi_* R\sset \k[\![x]\!].$

\item (The one-variable case, $n\!=\!1$)
 Let $\xi\!=\!a(x)\di_x,$ with $a(x)\in (x)^2\sset \k[\![x]\!].$ Denote $y(x)\!=\!e^\xi x.$
Then $\xi y(x)\!=\!a(y(x)),$ i.e. the power series $y(x)$ satisfies the ODE $y'=\frac{a(y)}{a(x)}.$
\\Indeed, $ \xi  y(x)\!=\!\xi e^\xi x\!=\!e^\xi \xi(x)\!=\!e^\xi a(x)\!=\!a(y(x)).$

\eee

\subsection{Algebraic power series}\label{Sec.Prelim.Algebraic.Series} Let $x=(x_1,\dots,x_n).$
Denote by $\k\bl x\br\sset \k[\![x]\!]$ the subring of algebraic power series over $\k[x].$ Namely, each $f(x)\in \k\bl x\br$ satisfies a polynomial
 equation $p(x,f(x))=0,$ for some $0\neq p(x,y)\in \k[x,y].$ Usually we take a minimal such equation, with $p(x,y)$ of the lowest $y$-degree.
  (In particular $p(x,y)$ is irreducible.)

\medskip

Occasionally we pass to the field of fractions.
\bel
$\Frac(\k\bl x\br)\!=\!\Frac(\k[x])\!\otimes _{\k[x]}\!\k\bl x\br.$
\eel
\bpr The part $\supseteq$ is immediate.
\\
For the part $\sseteq$ it is enough to show: $\frac{1}{f}\in \Frac(\k[x])\otimes _{\k[x]}\k\bl x\br$ for all
  $0\neq f\in \k\bl x\br.$ Take the minimal polynomial, $f^d+\cdots+a_0=0,$ where $0\neq a_0\in \Frac(\k[x]).$
   Thus $\frac{1}{f}=-\frac{f^{d-1}+\cdots+a_1}{a_0}\in \Frac(\k[x])\otimes _{\k[x]}\k\bl x\br.$
   \epr

\medskip

For $n=1$ this simplifies to $\Frac(\k\bl x\br)= \Frac(\k)\bl x\br[x^{-1}],$ i.e. Laurent polynomials, whose coefficients are algebraic power series.

For $\k$-local and Henselian, the ring  $\k\bl x\br$ is Henselian, closed under the derivations $\di_{x_1},\dots,\di_{x_n},$
 and admits compositions. For $\k$ a field this is well known.
 For $\k$ a local domain it is enough to observe: $\k\bl x\br=\Frac(\k)\bl x\br\cap \k[\![x]\!].$

\subsection{Differential extensions of $\pmb{\k[x]\sset R\sset  \k[\![x]\!]}$}\label{Sec.Prelim.Different.Extension}
 Let   $x\!=\!(x_1,\dots,x_n)$ and $R\sset\k[\![x]\!].$
\bee[\!\!\!\!\bf i.]
\item
 A power series  $f\in\k [\![x]\!]$ is called differentially definable over $R$ if the vector space
  $\Span_{\Frac(R)}(\{\di_{x_i}f\},\{\di^2_{x_ix_j}f\},\dots)$
  is finite dimensional over  $\Frac(R).$
  Equivalently, there exists a regular sequence of linear differential operators $D_1,\dots,D_n\in \Frac(R)[\{\di_{x_i}\}_i]$ satisfying:
   $D_1 f=0=\cdots=D_n f.$

 Denote by $D(R)$ the set of all the  power series differentially definable   over $R$.
  Then $D(R)\sseteq \k [\![x]\!]$ is a differential subalgebra. It is closed under derivations $\{\di_{x_i}\}$, antiderivatives $\int (\dots) dx_i,$
   and multiplicative inverses,
   \cite[Theorem 4]{PastorPillweinSinger2020}. The latter means: $D(R)\supseteq R_{(x)}.$
\\\mbox{Moreover, $D(R)\!\supseteq\! R^{hens}\!,$ \!the Henselization, see \cite[Theorem 5]{PastorPillweinSinger2020}.}

 \item
 In particular, $\k\bl x\br\sset D(\k[x])\sset \k[\![x]\!]$  is the subalgebra of $D$-finite power series (also called holonomic),
    \cite{StanleyDfinite}, \cite{KauersPaulebook}.
 The power series $e^x,sin(x),$ $cos(x),$ $ln(1+x),$ $arctan(x)$ are $D$-finite over $\k[x].$

We remark: $D(\k\bl x\br)=D(\k[x]).$

\item Similarly one defines $D^2(R):=D(D(R))$ and $D^{j+1}(R):=D(D^j(R)).$
 For example, $arcsin(x),e^{e^x-1},tan(x)\in D^2(\k[x])\smin D(\k[x]).$

The  rings $D^\bullet(R)$ are in general neither  local, nor Henselian.
However, their $(x)$-localization  and $(x)$-Henselization are controllable, $D(R)_{(x)}\sseteq D(R)^{hens}\sseteq D^2(R).$
   (Apply \cite[Proposition 6]{PastorPillweinSinger2020} to the ring $D(R)$.)

\item The rings $D^\bullet(\k[x])$ do not admit compositions. But a weaker property holds:
\bei
\item The rings $D^\bullet(\k[x])$ admit compositions by the elements of $\k\bl x\br,$ \cite[Proposition 6]{PastorPillweinSinger2020}.
\item If  $f\in D^i(\k[x])$ and $g\in (x)\cdot D^j(\k[x]),$ then $f\circ g\in D^{i+j}(\k[x]),$  \cite[Theorem 10]{PastorPillweinSinger2020}.
\eei
\item
The natural (inductive) limit
 $D^{\infty}(\k [x])\!:=\!\bigcup_{j\geq 0} D^j (\k [x])\sset \k[\![x]\!] $
 is a differential subalgebra that admits compositions, is
  local,  and contains all the   antiderivatives, \cite[page 16]{PastorPillweinSinger2020}.
 The ring  $D^{\infty}(\k [x])$ is regular, local, Henselian. This holds because $(D^j (\k [x]))^{hens}\sseteq D^{j+1} (\k [x]),$ see \S\ref{Sec.Prelim.Different.Extension}.i.

Another name for $D^{\infty}(\k [x])$ is the Picard-Vessiot closure of $\k[x]$ inside $\k[\![x]\!].$

\item
A power series $f\in\k[\![x]\!]$ is called differentially algebraic   if the field
  $\k(x,\{\di_{x_i}f\},\{\di^2_{x_ix_j}f\},\dots)$ has finite transcendence degree over $\k$.
 Thus $f$ satisfies a polynomial ODE $P(x,y,y',\dots,y^{(N)})=0,$ and the minimal such  $N$ is called {\em the differential order of $f.$}

  The set of differentially algebraic power series forms a differential   subalgebra
  $ D^{alg}(\k[x])\sset  \k[\![x]\!].$
   All $D^{\infty}$-power series are differentially algebraic, \cite[Corollary 29]{PastorPillweinSinger2020},
  i.e. $D^{\infty}(\k[x])\sset D^{alg}(\k[x]).$
   The algebra $D^{alg}(\k[x])$ is local and admits compositions. Moreover, the set of algebras $D^{alg}(\k[x])$ (for $n=1,2\dots$)
     is a Weierstrass system over $\k.$ In particular it admits the implicit function theorem, \cite[(5.2) Theorem]{vandenDries1988}.
  In particular it is a Henselian   ring.
\eee

\subsection{The  Lefschetz principle for the passage from $\k$ to $\C$}\label{Sec.Prelim.Lefschetz}
 Let $x=(x_1,\dots,x_n).$
Fix a derivation $\xi=\sum a_i \di_{x_i}\! \in\! \Der_\k^{nilp}(\k[\![x]\!])$ and let $\{c_\bullet\}\sset \k$
  the set of all the coefficients of all $\{a_i\}.$
 Thus $e^\xi x_i\!\in\! \Q(\{c_\bullet\})[\![x]\!]$ for each $i.$

 The field extension
   $\Q\sseteq \Q(\{c_\bullet\})$ is (at most) countably generated. Therefore there exists an embedding
 of fields $\phi: \Q(\{c_\bullet\})\into \C.$
   Accordingly we get a complex derivation
  $\xi^\C:=\phi(\xi)\in \Der_\C^{nilp}(\C[\![x]\!]).$ Then $\phi(e^\xi x)\in \Q(\phi\{c_\bullet\})[\![x]\!]\sset \C[\![x]\!].$
 Hence   the power series $y_\C(x):=\phi(e^\xi x)$ satisfy some implicit-function
  equations\footnote{Here $F_\C\in \Q(\phi\{c_\bullet\})[ x,y  ]$ or in $\Q(\phi\{c_\bullet\})\bl x,y  \br$
  or in other subrings of $\Q(\phi\{c_\bullet\})[\![ x,y   ]\!].$}  $F_\C(x,y )=0 $   \iff
   the initial power series $e^\xi x$ satisfy the corresponding equations over $\k.$
Finally, if $0=F(x,y_\C(x))$ for some $0\neq F(x,y)\in \C\{x,y\}$, then $y_\C(x)$ satisfies a corresponding (polynomial/algebraic/analytic)
  equation over $ \Q(\phi\{c_\bullet\})$ as well.
 (E.g. by taking a Hamel basis of $\C$ over $\Q(\phi\{c_\bullet\})$.)

\medskip

In this way we conclude: if $\phi(e^\xi x)\in \C\bl x\br$ then $e^\xi x\in \k\bl x\br.$
 Similarly, if $y_\C$ satisfies a (non)linear ODE over $\C[\![x]\!]$ then so does $e^\xi x $ over   $ \Q(\phi\{c_\bullet\})[\![x]\!].$
 Hence, whenever we work with $\k[x],$  $\k\bl x\br,$   $D^\bullet(\k[x]),$ $D^\infty(\k[x]),$  $D^{alg}(\k[x]),$ the reader  (if needed)
  can assume $\Q\sseteq\k\sseteq\C.$

\subsection{The action $e^\xi\circlearrowright R$ vs the solvability of the ODE $\pmb{x'=\xi(x)}$}\label{Sec.3.Exponentiability.as.ODE.solvability}
 Let  $x=(x_1,\dots, x_n).$ Take a differential $\k[t]$-subalgebra   $\k[t,x]\sset R\sseteq \k[t][\![x]\!],$ not necessarily Noetherian.      Denote $R_o:=R\cap \k[\![x]\!].$
  For each $t_o\in \k$ define the evaluation map $ev_{t_o}:\k[t][\![x]\!]\to \k[\![x]\!] $ by $f(t,x)\to f(t_o,x).$  Thus $ev_{t_o}(R)\supseteq R_o.$

In this subsection   we assume the following properties.
 \bee[\bf a.]
 \item ($R$ is closed under $t$-evaluations) $ev_{t_o}(R)=R_o$ for every $t_o\in \k.$
 \item ($R$ is ``saturated" inside $\k[t][\![x]\!]$) \quad  $R[\frac{1}{f}]\cap \k[t][\![x]\!]=R$ for each $0\neq f\in R.$

 \eee
\bex\label{Ex.rings.for.section.2.5}
\bee[\bf i.]
\item (Negative examples) The ring $R=\k[t,x]+(t)\cdot \k[t][\![x]\!]$ violates a. and b.
 The ring  $R=\k[t,x]+e^{2x}\cdot \k[t,x,e^{ x}]$ satisfies a. but violates b.

\item Assumption b. implies: $R$ is ``x-local", i.e. all elements of the set $\{1\}+(x)$ are invertible in $R.$
  Thus the rings  $D^j_x(\k[t,x]),$ for $1\le j<\infty,$ do not satisfy b.
 \item (Positive examples)
The rings $\k[t]\bl x\br,$ \ $\k[t]\{x\}:=\Frac(\k)\{t,x\}\cap \k[t][\![x]\!],$ \ $D^\infty_x(\k[t,x]),$ $D^alg_x(\k[t,x]),$
 satisfy these assumptions.
\eee
\eex

The evaluation $R\stackrel{ev_{t_o}}{\to}R_o$ defines the evaluation of derivations,
  $\Der_{\k[t]}(R)\stackrel{ev_{t_o}}{\to}\Der_\k( \k[\![x]\!]),$ by $\xi=\sum a_i(t,x)\di_{x_i}\to \xi|_{t_o}:=\sum a_i(t_o,x)\di_{x_i}.$

\bel
The image of this map consists of $R_o$-derivations,
 $\Der_{\k[t]}(R)\stackrel{ev_{t_o} }{\twoheadrightarrow} \Der_{\k }(R_o).$
 \eel
\bpr
Let $\xi=\sum_i a_i(x,t)\di_{x_i}\in \Der_{\k[t]}(R).$ Then $a_i(x,t)=\xi(x_i)\in R.$ Therefore $a_i(x,t_o) \in R_o.$
 Hence $ev_{t_o}(\xi)\in \Der_{\k }(R_o).$

To verify the surjectivity, take any $\xi\in \Der_{\k }(R_o),$ and extend it $t$-linearly to  $\xi\in \Der_{\k[t]}( \k[t][\![x]\!]).$
 We claim: $\xi\in \Der_{\k[t]}(R).$ Indeed, present $\xi=\sum_i a_i(x)\di_{x_i},$ where $a_i(x)\in R_o\sset R.$ Therefore
  $\xi(R)\sseteq R_o\cdot (\di_1\dots \di_n)(R)\sseteq R.$
   \epr

\medskip

Take  $\xi\in \Der^{nilp}_{\k}(R_o)$ and extend it to  $\xi\in \Der^{nilp}_{\k[t]}(R).$
\bprop\label{Thm.TFAE}  The following properties are equivalent.
\bee
\item $e^\xi x\in R_o.$
\item $e^{t_o\cdot\xi}x\in R_o$  for each $t_o\in \k.$

\item The ODE $x'\!=\!\xi(x)$ is solvable over $R,$ i.e. $x(t)\!\in\! R$ for the unique solution $x(t)\!=\!e^{t\xi}x.$

\item (Assuming $R_o,R$ admit compositions, \S\ref{Sec.Prelim.Notations.Conventions}) \ $e^{\xi}$ acts on $R_o,$  $e^{t\xi}$ acts on $R.$
\eee
\eprop
\bpr
The equivalence $3 \Leftrightarrow4$  is trivial.

The implication $3\Rrightarrow 2$ follows by the evaluation map, $ev_{t_o}(e^{t\xi})=e^{t_o\xi}.$

The implication $2\Rrightarrow 1$  is trivial.

We prove $1\Rrightarrow 3.$
 Take $\xi\in  \Der^{nilp}_{\k}(R_o)$ and the corresponding  power series $e^{t\xi}x\in \k[t][\![x]\!].$
 Consider the quotient $M:=\quots{R+R\cdot e^{t\xi}x }{R}\sset \quots{\k[t][\![x]\!]}{R}.$
 This is an $R$-submodule, it is finitely generated by the power series $e^{t\xi}x_1,\dots,e^{t\xi}x_n.$

 We claim: the fibre $M\otimes_{\k[t]}\quots{\k[t]}{(t-1)}$ vanishes. Indeed, one has:
\beq
M\ni [e^{t\xi}x]\stackrel{e^\xi x\in R}{=}[e^{t\xi}x-e^{\xi}x]=[(t-1)\sum^\infty_{j=1}c_j(t,x)],
\eeq
 for some polynomials $c_j(t,x)\in (x)^j\sseteq \k[t,x].$
 \bei
 \item If $\sum^\infty_{j=1}c_j(t,x)\!\in\! M$ then immediately:
 \beq
     [(t-1)\sum^\infty_{j=1}c_j(t,x)]\!=\!(t-1)\cdot [\sum^\infty_{j=1}c_j(t,x)]\!=\!0\!\in\! M\!\otimes\!\quots{\k[t]}{(t-1)}.
     \eeq
 \item Without this assumption one argues inductively:
 \beq
 M\!\otimes\!\quot{\k[t]}{(t-1)}\ni [e^{t\xi}x] \stackrel{\forall N}{=}0+[(t-1)\sum^\infty_{j=N}c_j(t,x)].
 \quad\quad \text{Thus\ } [e^{t\xi}x]\in \cap^\infty_{j=1}[ (x)^j] =0.
 \eeq
 \eei
As $M$ is finitely generated,  the support $\Supp(M)$ is a closed subscheme of $\Spec(R),$  and $\Supp(M)\ssetneq \Spec(R).$
  Thus $\Supp(M)\sseteq V(f)\sset \Spec(R),$ for some  $0\neq f\in R.$
  Passing to the ring of fractions $R[\frac{1}{f}]$ we get the vanishing: $0=M[\frac{1}{f}]\in mod$-$R[\frac{1}{f}].$
\\
Therefore $e^{t\xi}x\!\in\! R[\frac{1}{f}]\!\cap\! \k[t][\![x]\!].$
   Finally  $R[\frac{1}{f}]\!\cap\! \k[t][\![x]\!]\!=\!R$ by the saturatedness assumption.
 \epr
\bex
\bee[\bf i.]
\item We restate this lemma for the ring $\k[t]\bl x\br$:
 $e^\xi $ acts on $\k\bl x\br$ \  \iff  \  $e^{t\xi} $ acts on $\k[t]\bl x\br$  \  \iff  \ the ODE $x'=\xi(x)$ is solvable over $\k[t]\bl x\br.$
\\
And similarly for the rings of Example \ref{Ex.rings.for.section.2.5},  $S\{x\},$ $D^\infty(S[x]),$ $D^{alg}( S[x]),$ where $S=\k[t].$

\item Part 2 of Proposition \ref{Thm.TFAE} implies:  the subset of exponentiable derivations is a cone inside the $\k$-module  $\Der^{nilp}_\k(R_0).$

    However, this cone is not a submodule. E.g. let $n=1$ and $R=\k\bl x\br$ and $\xi_i=\frac{x^3}{1+a_i x}\frac{d}{d x}$
    for some constants $a_i\in \k.$ Thus $\xi_1,\xi_2$ are exponentiable derivations, $e^{\xi_1}x,e^{\xi_2}x\in\k\bl x\br, $
     see Corollary \ref{Thm.exp.for.exact.form}.      But
     \beq
     \xi_1+\xi_2=\frac{x^3(2+a_1 x+a_2 x)}{(1+a_1 x)(1+a_2 x)} \cdot \frac{d}{d x}=\frac{2x^3}{1+\frac{a_1+a_2}{2}x-(\frac{a_1-a_2}{2})^2x^2+x^3(\dots)}\cdot\frac{d}{d x}.
     \eeq
      And this is non-exponentiable for $a_1\neq a_2,$ again by   Corollary \ref{Thm.exp.for.exact.form}.
\eee
\eex

\section{Algebraicity/transcendence in the one-variable case, $n=1$} 

\subsection{\!Computational examples  for $\k$   a field}\label{Sec.Exp.Der.n=1}\label{Sec.n=1.Computing.exp}

 Take a  nilpotent derivation $\xi\!\in\! (x)^2\cdot \Der_\k(\k[\![x]\!]) $ of order $p+1\ge 2.$  We pass to differential forms
via the condition $\om(\xi)=1.$ Thus we get $\om=\ta(x)\cdot dx,$
where $\ta(x)\in  \k[\![x]\!] $ has $x$-order $(-p-1)\le -2.$
 Take the residue at the origin, $res_o\om\in  \k .$ Thus we can present $\om=[\di_x \frac{1}{c\cdot  a(x)^{p} }+\frac{res_o\om}{x}]dx,$
  where $0\neq c\in \k$ and $a(x)-x\in (x)^2\sset  \k [\![x]\!].$

\bel\label{Thm.Exp.Der.Case.n=1}
\bee
\item
The power series $y(x)=e^\xi x\in \k[\![x]\!]$  is determined by the    equation
\beq\label{Eq.n=1.Exp.Der}
  \frac{1}{c\cdot a(y)^{p}} - \frac{1}{c\cdot a(x)^{p}}  +res_o\om\cdot ln\frac{y}{x}=1.
  \eeq

 \item  In particular, if $res_o\om\!=\!0$  then  $e^\xi x\!=\!a^{-1}[\frac{a(x)}{\big(1+c\cdot a(x)^{p}   \big)^\frac{1}{p}}].$
\eee
\eel
 Here $a^{-1}[\dots]$ is the inverse of the power series $a(x)$.
\bpr
Consider the ODE $x'=\xi(x).$ Its solution is $e^{t\xi}x,$ thus $e^\xi x$ is the time-one flow. Passing to the Pfaffian form we get
 $1=\int^y_x\om=\frac{1}{c\cdot a(z)^p}\vert^y_x+res_o\om \cdot ln(z)\vert ^y_x.$ Hence the first statement.

 The second statement is immediate.
\epr

\bex\label{Ex.n=1.Computations}
\bee[\bf i.]
\item
Let   $\xi=\frac{c}{p}\cdot x^{p+1} \di_{x}$  for $b\in \k$ and $p\ge1.$ Thus   $\om=-\frac{\di_{x} (x^{-p})}{c} dx .$
  Lemma \ref{Thm.Exp.Der.Case.n=1} gives:
$e^\xi x = \quots{x}{\big(1-c \cdot x^p\big)^\frac{1}{p}}.$

By the multiplicativity of $e^\xi $ we get for any power series: $e^\xi  f(x) \!=\!f(\quots{x}{\big(1-c\cdot x^p\big)^\frac{1}{p}} ).$

In particular, the operator   $e^\xi$ acts on $\k\langle x\rangle,$   on $ D^j(\k[x])$ for $j\le\infty,$ and on
 $D^{alg}(\k[x]).$

\item Let $\xi=c\cdot sin^2(x)\cdot \di_x,$ thus $\om=-\frac{1}{c} \cdot\di_x \frac{1}{tan(x)}\cdot  dx.$
 We get: $e^\xi x=arctan\big[\frac{tan(x)}{1-c\cdot tan(x)}\big].$

 In this case $e^\xi x\in D^2(\k[x]).$ Indeed, $(e^\xi x)'=\frac{1}{(1+\frac{c^2}{2})-c\cdot sin(2x)-\frac{c^2}{2} cos(2x)}.$
\item
Let $\xi=\frac{c}{p}\cdot\frac{sin^{p+1}(x)}{cos(x)}\cdot \di_x,$ thus $\om=-\frac{1}{c} \cdot\di_x \frac{1}{sin^p(x)}\cdot  dx.$
 We get: $e^\xi x=arcsin\big[\frac{sin(x)}{\sqrt[p]{1-c\cdot sin^p x}}\big].$

In this case $e^\xi x\in D^2(\k[x]).$ Indeed, denote $f(sin(x)):= e^\xi x,$ i.e. $f(t)=arcsin\big[\frac{t}{\sqrt[p]{1-c\cdot t^p}}\big].$
 Then $f'(t)\in \k\bl t\br,$ hence $f (t)\in  D(\k[t]).$

\item Let $\xi=\frac{c}{p}\cdot\frac{(e^x-1)^{p+1}}{  e^x}\cdot \di_x,$ thus $\om =-\frac{1}{c} \cdot \di_x \frac{1}{(e^x-1)^p}\cdot \di_x.$
We get:  $e^\xi x=ln\big[1+\frac{e^x-1}{\sqrt[p]{1-c\cdot (e^x-1)^p}}\big].$

In this case $e^\xi x\in D^2(\k[x]).$ (The proof is similar.)
 \eee\eex
\bex\label{Ex.Lambert.Function}
Despite these ``positive"  cases, in general $e^\xi x$ need not belong even to $D^\infty(\k[x]).$
  As an example, let $\xi\!=\!-\frac{x^2e^{-x}}{x-1}\di_x,$ here $\frac{x^2e^{-x}}{x-1}\!\in \! D(\k[x]) .$
 {\bf Claim:} $e^\xi x\!\not\in\! D^\infty(\k[x]).$
\bpr
  Here $\om\!=\!\di_x (\frac{e^{x} }{x  })dx.$
   Thus equation \eqref{Eq.n=1.Exp.Der} is: $\frac{e^{ y}}{y  }-\frac{e^{ x}}{x  }\!=\!1.$
    Therefore $y(x)$ satisfies the equation     $y\cdot e^{-y}\!=\!\frac{x\cdot e^{-x}}{1+x\cdot e^{-x}}.$
     But the  equation $z\cdot e^z=t$ defines the Lambert function $z(t)\!=\!W(t),$
     which is known to be not in $D^\infty(\k[t]).$ (Theorem 1 in \cite{Bostan-Jimenez-Pastor} or \cite{Bronstein.Corless.Davenport.Jeffrey})
      Finally, we can present: $y(x)\!=\!-W(-\frac{x\cdot e^{-x}}{1+x\cdot e^{-x}}).$
       Here $\frac{x\cdot e^{-x}}{1+x\cdot e^{-x}}\!\in\! D^2(\k[x])$ and the ring $D^\infty(\k[t])$ admits compositions.
        Therefore $W(-\frac{x\cdot e^{-x}}{1+x\cdot e^{-x}})\!\not\in D^\infty(\k[x]).$
\epr
   \eex

\beR
  When $\k$ is not a field, but a local domain over $\Q,$ the same formulae hold. But now the results are in the ring $\Frac(\k)[\![x]\!].$
\eeR

\subsection{Several corollaries}

\bcor\label{Thm.no.WPT.for.Dinfty}
The ring $D^\infty(\k[x,y]),$ being regular, local and Henselian, does not admit the Weierstrass preparation.
\ecor
\bpr
Take a derivation $\xi\in x^2\cdot \Der_\k(D^\infty(\k[x])).$ Denote $y(x)=e^\xi x=x\cdot (1+\ty(x)).$ Equation \eqref{Eq.n=1.Exp.Der} becomes now:
\beq
\frac{1}{c\cdot a(x\cdot (1+\ty(x)))^{p}} - \frac{1}{c\cdot a(x)^{p}}  +res_o\om\cdot ln  (1+\ty(x)) =1 .
\eeq

Multiplying by $c\cdot a(x)^{p}$ we get:
\beq
 0= F(x,\ty(x)):=\frac{  a(x)^{p}}{ a(x\cdot (1+\ty(x)))^{p}} - 1  +c\cdot a(x)^{p}\cdot \big(res_o\om\cdot ln  (1+\ty(x)) -1\big).
\eeq
Note that $F(x,\ty)\in D^\infty(\k[x,y]).$ Observe: $\di_\ty F|_{o,o}\neq0.$

Suppose the ring  $D^\infty(\k[x,y]) $   admits the Weierstrass preparation. Then we can present $ F(x,\ty)=u(x,\ty)\cdot (\ty-h(x)),$
 for some $h(x)\in D^\infty(\k[x])$ and $u(o,o)\neq0.$ But then we get: $e^\xi x\in D^\infty(\k[x])$ for any $a(x)\in D^\infty(\k[x]).$
  In contradiction to Example \ref{Ex.Lambert.Function}.
\epr

\medskip

Unlike the negative Example \ref{Ex.Lambert.Function}, for differentially-algebraic power series, $R:=D^{alg}(\k[x]),$ the result is positive.
\bcor\label{Thm.n=1.exp.for.Diff.algebraic}
\bee
\item
$\Der^{exp}_\k(R)=(x)^2\cdot \Der_\k(R),$ i.e. $e^\xi x\in R$ for every $\xi\in (x)^2\cdot \Der_\k(R).$
\item For algebraic derivations, $\xi\in (x)^2\cdot \Der_\k(\k\bl x\br),$ the differential order of $e^\xi x$ is 1.
\eee
\ecor
\bpr
\bee
\item
Proceed as in the proof of Corollary \ref{Thm.no.WPT.for.Dinfty} to get $F(x,\ty)\in R.$
 And observe that $R$ admits the implicit function theorem, \S\ref{Sec.Prelim.Different.Extension}.vi.
 \item
 Present $\xi=b(x)\cdot \di_x,$ then the power series $y(x)=e^\xi x$ satisfies  $y'=\frac{b(y)}{b(x)},$ by \S\ref{Sec.Prelim.Derivations}.iv.
 Consider the element  $w:=b(x)\cdot z-b(y)\in \k\bl x,y,z\br.$ Take its minimal polynomial, $\sum P_j(x,y,z) w^j=0.$
  Here $P_j\in \k[x,y,z]$ and $P_0\neq0$ by minimality. Therefore $y(x)$ satisfies the (nontrivial) polynomial ODE $P_0(x,y,y')=0.$
\epr
\eee

 Start from  a local domain  $\k\supseteq\Q$ and pass to the field of fractions $\Frac(\k).$
 For algebraic power series,  $R\!=\!\k\bl x\br,$
 it is easy to describe the set of exponentiable derivations. Present  $\xi\!\in\! (x)^2\cdot \Der_\k(\k\bl x\br)$
 as in Lemma \ref{Thm.Exp.Der.Case.n=1}. We use those notations.
 \bcor\label{Thm.exp.for.exact.form}
 \bee
 \item  $\xi$ is exponentiable \iff $res_o\om=0$ and $a(x)^p\!\in \!\Frac(\k)\bl x\br.$
\item
 The $(x)$-adic closure of the set of exponentiable derivations, $\overline{\Der^{exp}_\k(\k\bl x\br)}\sset (x)^2\cdot \Der_\k(\k[\![ x]\!]),$
is defined  by the single condition $res_o\om=0.$
\eee
\ecor
In particular the exponentiability of $\xi$ is not finitely determined (due to the condition  $a(x)^p\!\in \!\Frac(\k)\bl x\br$).
 For any $\xi\in  \Der^{exp}_\k(\k\bl x\br)$ and any $N\in \N$
 there exists $\de\in (x)^N\cdot \Der_\k (\k[ x])$ such that $\xi+\de\not\in \Der^{exp}_\k(\k\bl x\br).$

Part 1 can be restated: $\xi$ is exponentiable \iff the corresponding differential form $\om$ is exact (over $R$). Thus the set of exponentiable
 differential forms is a $\k$-module.

\subsection{The collapse from $e^\xi x\in D^\infty(R)$ to $e^\xi x\in  R^{hens}$}

\subsubsection{\!\!A special implicit function equation for   $e^\xi x\!\in\! D(R)$}
\mbox{Let $\xi \!=\!a(x)\di_x,$ with $a(x)\!\in\! R\sseteq\k[\![x]\!].$}
\bel\label{Thm.exp.in.Der(R).vs.R^hens}
 If $e^\xi x\in D(R)$ then $y(x)=e^\xi x$ satisfies an equation $F(x,y)=0,$ for a power series $0\neq F(x,y)\in R_x\otimes_\k R_y\sset \k[\![x,y]\!].$
\eel
\bpr Denote $y(x)=e^\xi x.$
 If $y(x)\!\in\! D(R)$ then it satisfies a linear ODE,
 $\sum^r_{j=0}c_j(x) y(x)^{(j)}=0,$ for some $c_j(x)\!\in\! R.$
  Using $\xi_x y(x) =a(y(x) ),$ see \S\ref{Sec.Prelim.Derivations}.iv., and its higher versions, $\xi^j_x y(x) =\big(a(y) \frac{d}{dy}\big)^jy|_{y(x)}, $ we get
   $F(x,y(x))=0$ for some
 \beq
 F(x,y):=\sum^r_{j=0} \tc_j(x)\cdot \big(a(y) \frac{d}{dy}\big)^jy\in R_x\otimes_\k R_y \sset \k[\![x,y]\!].
 \eeq
  Finally, $F(x,y)\!\neq\!0,$ because $\ord_y (a(y) \frac{d}{dy})^jy\!=\!j\!\cdot\! (\ord_y a(y)-1)\!+\!1.$ (Note that $a(y)\!\in\! (y)^2.$)
\epr

\subsubsection{} Let $\xi =a(x)\di_x$ with $a(x)\in (x)^2\sset \k\bl x\br.$
For  the simple cases of Example \ref{Ex.n=1.Computations} we saw: either $e^\xi x\in \k\bl x\br$ or $e^\xi x\not\in D(\k[x]).$
 A much stronger property holds, with $R^{hens}$-the Henselization.

\bthe\label{Thm.D^infty.to.Rhens}
If $e^\xi x \in D^\infty(R),$ then $e^\xi x \in R^{hens}.$
\ethe
The converse  is trivial, as $D(R)\supseteq R^{hens},$ see \S\ref{Sec.Prelim.Different.Extension}.i.

This theorem does not hold in the $D$-finite case,   $a(x)\in D(\k[x]),$ in view of Example \ref{Ex.Lambert.Function}.
\bpr First we observe (for any subring $S\sset\k[\![x]\!]$): if $e^\xi x\in D(S)$ then $e^\xi x\in S^{hens}.$
Indeed, by Lemma \ref{Thm.exp.in.Der(R).vs.R^hens} the power series $y(x):=e^\xi x$ satisfies $F(x,y(x))=0$ for some
  $0\neq F(x,y)\in S_x\otimes_\k\k\bl y\br\sseteq S^{hens}_x\bl y\br.$
 Invoking the Weierstrass preparation theorem\footnote{Here $S$ is not necessarily Noetherian. See \cite{Moret-Bailly} for this case, and for the
    more classical references.} for the Henselian ring $S^{hens}_x\bl y\br,$ we can present $F(x,y)=u(x,y)\cdot p(x,y),$
     where $0\neq p(x,y)\in S^{hens}_x[y],$
     while $u(x,y)|_{(o,o)}\neq0.$ Therefore $p(x,y(x))=0,$ hence $y(x) \in S^{hens}.$

\medskip

 It is enough to prove: if $y(x)\in D^l(R) $ for some $l\ge2,$
 then $y(x)\in  D^{l-2}(R)^{hens}.$  Below we use the field of fractions $\K_l:=\Frac(D^l(R)).$
\bee[\bf Step 1.]
\item Let $y(x)\in D^l(R),$ then $y(x)\in D^{l-1}(R)^{hens},$ by the observation above.
 Namely, $P(y(x))=0$ for a polynomial $0\neq P(y)\in \K_{l-1}[y].$

W.l.o.g. we take $P(y)$ as the minimal monic polynomial of $y(x).$ Thus $P(y)$ is irreducible, and all its roots, $\{y_j(x)\}\sset \overline{\K_{l-1}},$
 are distinct. One of them is the power series $y(x),$ and it satisfies the ODE $y'=a(y)/a(x),$ see \S\ref{Sec.Prelim.Derivations}.iv.
  We claim: all the roots satisfy this ODE.

Indeed, apply $\frac{d}{dx}$ to the identity $P(y(x))=0,$ to get $Q(y(x))=0$ for a polynomial
  $Q(y):=a(x)\cdot \di_x P(y)+a(y)\cdot \di_y P(y)\in {\K_{l-1}}[y].$ Here $Q(y)\neq0,$ e.g. because
 \beq
 \ord_y (a(y)\di_y P(y))\ge \ord_y a(y)+\ord_y P(y)-1>\ord_y P(y).
 \eeq
But $P(y)$ is the minimal polynomial of $y(x).$ Therefore $P(y)$ divides $Q(y).$ In particular, $Q(y_j(x))=0$ for all the roots $y_j(x).$

On the other hand, apply $\frac{d}{dx}$ to the identity $P(y_j(x))=0$ to get: $\di_x P(y)|_{y=y_j(x)}+\di_y P(y)|_{y=y_j(x)}\cdot y_j(x)'=0.$
 Together with $Q(y_j(x))=0$ we get:
 \beq
 \di_y P(y)|_{y=y_j(x)}\cdot \Big(y_j(x)'-\frac{a(y_j(x)}{a(x)}\Big)=0.
 \eeq
Finally, $ \di_y P(y)|_{y=y_j(x)}\neq0$ (separability), as the roots are simple and $char(\k)=0.$
 Therefore $y_j(x)'=a(y_j(x))/a(x).$

\medskip

\item Expand $P(y)=\sum_k b_k\cdot y^k\in \K_{l-1}[y].$
 As the polynomial $P(y)$ is monic, its coefficients $b_\bullet$ are the elementary symmetric polynomials in the roots $\{y_j(x)\}.$
 These symmetric polynomials generate the subring $\k[\uy]^{sym}\sset \k[\uy].$ We take another system of generators of $\k[\uy]^{sym}:$
 \beq
 \si_k(\uy):=\sum_j y^k_j,\quad \text{for } k=1,\dots, \deg P(y).
 \eeq
 Then $\si_k(\uy)\in\k[b_\bullet]\sset  \K_{l-1},$ i.e. these generators satisfy linear ODE's over $\K_{l-2}.$
 Present these ODE's in the form:
\beq
\{ \sum_i  c_{k,i}\cdot (a(x)\di_x )^i \si_k(\uy)=0\}_k,\quad \text{for some }\quad  c_{k,i}\in \K_{l-2}.
\eeq
Each root $y_j(x)$ satisfies the ODE $y'=a(y)/a(x),$ therefore
\beq
(a(x)\di_x )^i \si_k(\uy)=(\sum_j a(y_j)\di_{y_j})^i\si_k(\uy)= \sum_j (a(y_j)\di_{y_j})^i y_j^k.
\eeq
Denote this symmetric power series  by $q_{k,i}(\uy)\in \k\bl \uy\br$  Altogether, the roots $\uy(x)$ satisfy
  the system of    equations
 \beq\label{Eq.Equations.to.prove.Regularity}
\{\sum_i c_{k,i}\cdot q_{k,i}(\uy)=0\}_k,\qquad q_{k,i}(\uy).
 \eeq
 Denote $n=\deg_y P(y),$ thus we have $n$ power series in $n$ variables.

\

\hspace{-1.7cm}{
\parbox{16cm}{It remains to prove: \eqref{Eq.Equations.to.prove.Regularity} is a regular sequence   in the ring $\K_{l-2}\bl y\br.$
 Then we will get $y_j(x)\in \overline{\K_{l-2}}$ for each $j.$ In particular $y(x)\in D^{l-2}(R)^{hens}.$
  Iterating Steps 1 and 2 we get eventually:  $y(x)\in R^{hens}.$ Which is exactly the statement.}}

\medskip

\item  We prove:  equations \eqref{Eq.Equations.to.prove.Regularity}  form a regular sequence in the ring $\K_{l-2}\bl y\br.$

\bei
\item {\em The naive approach.}
It is enough to verify: their leading terms (of lowest $\uy$-orders)
 form a regular sequence. (See e.g. Proposition 15.15 in \cite{Eisenbud}.) We have:
 \beq\label{Eq.leading.terms}
 l.t.(q_{k,i}(\uy))=\sum_j(l.t.[ a(y_j)]\cdot \di_{y_j})^i y^k_j=c\cdot \sum_j y^{N_{k,i}}_j.
 \eeq
 Here $c=const\neq0$ and $N_{k,i}=(\ord_y a(y)-1)i+k.$
 Therefore the leading term of $k$'th equation in \eqref{Eq.Equations.to.prove.Regularity} (of lowest $\uy$-order) is
 $ c\cdot c_{k,i_{k,min}}\cdot \sum_j y^{N_{k,i_{k,min}}}_j.$ Here $i_{k,min}$ is the lowest index for which $c_{k,i}\neq0.$

We would like to deduce:  the polynomials $\{\sum_j y^{N_{k,i_{k,min}}}_j\}_k$ form a regular sequence.
  It is enough to verify (for an algebraically closed field $\k=\bk$) that their zero-locus is of dimension zero.
   And for that it is enough to verify: the Jacobian
   matrix of these polynomials is non-degenerate on their zero locus.

 Suppose all the indices $\{i_{k,min}\}_k$  coincide. Then $N_{1,i_{1,min}},N_{2,i_{2,min}},\dots,N_{n,i_{n,min}}$ is a
 sequence of consecutive integers.
   Hence the Jacobian matrix is the Vandermonde matrix, and its rank drops only at the points where $y_j=y_{j'}$ or $y_j=0.$
 Thus the leading terms form a regular sequence.

But the indices  $\{i_{k,min}\}_k$ need not coincide. And in this case the   polynomials $\{\sum_j y^{N_{k,i_{k,min}}}_j\}_k$
 might not form a regular sequence. (E.g. take $n=2$ and
 the polynomials $y_1+y_2$, $y^3_1+y^3_2$.)

\medskip

 Below we will adjust the system \eqref{Eq.Equations.to.prove.Regularity} to achieve: all $\{i_{k,min}\}_k$ coincide.

\item{\em The corrected approach.}
 The coefficients $\{c_{k,i}(x)\}$ belong to the field $\K_{l-2}.$ By scaling each equation we can assume: all $\{c_{k,i}(x)\}$ are regular at $x\!=\!0,$
  and for each $k$ at least one $c_{k,i}(0)$ is non-zero.
 Therefore equations \eqref{Eq.Equations.to.prove.Regularity} can be $x$-expanded:
\beq\label{Eq.Equations.to.prove.Regularity.split}
  \sum_{i\ge i_{k,min}} c_{k,i}(0)\cdot q_{k,i}(\uy)+x\cdot(\cdots)=0,\qquad c_{k,i_{k,min}}(0)\neq0.
\eeq
It is  enough to verify: the evaluations of  \eqref{Eq.Equations.to.prove.Regularity.split}  at $x\!=\!0$ form a regular sequence in $\k\bl\uy\br.$
\bei
\item
 If all the indices $\{i_{k,min}\}_k$ coincide, then we invoke the naive argument (above), with the Jacobian matrix and Vandermonde determinant.
\item
 Suppose these indices do not coincide. Apply $a(x)\di_x$ to the equation with the smallest $i_{k,min}.$
  The splitting of  \eqref{Eq.Equations.to.prove.Regularity.split}  is preserved, because $\ord_x a(x)\ge2.$ The new equation is:
\beq
 \sum_{i\ge i_{k,min}} c_{k,i}(0)\cdot  q_{k,i+1}(\uy)+x\cdot(\cdots)=0.
\eeq

And the power series $\uy(x)$ satisfy this new system of equations.

Iterate this procedure (if needed) to get equations of type \eqref{Eq.Equations.to.prove.Regularity.split}, where all the indices  $\{i_{k,min}\}_k$
 coincide. And now we invoke the naive argument, with the substitution $x=0.$
\epr
\eei

\eei
\eee

\bcor\label{Thm.If.D^infty.then.algebraic}
For $\xi\in (x)^2\cdot\Der_\k(\k\bl x\br)$   we have: if $e^\xi x \in D^\infty(\k[x])$ then $e^\xi x \in \k\bl x\br.$
\ecor

\subsection{$D^\infty$-classification for polynomial derivations}
\mbox{Let   $\xi\!=\!a(x)\!\cdot\! \di_x$ with $a(x)\!\in\! (x)^2\sset \k[x].$}

\bthe\label{Thm.exp.Der.is.D1.iff.xi=x^p}
 $e^\xi x\in D^\infty(\k[x])$ \iff $\xi =c\cdot x^p\cdot \di_x,$ for some $c\in \k.$
\ethe
Thus Example \ref{Ex.n=1.Computations}.i is the only occasion of $D^\infty$-result.

 This theorem cannot be extended to the case of algebraic power series, $a(x)\in \k\bl x\br,$
 due to part 2 of Corollary \ref{Thm.exp.for.exact.form}.
\bpr
The part $\Lleftarrow$ is Example \ref{Ex.n=1.Computations}.i. We prove the part $\Rrightarrow$.
We assume $a(x)\neq c\cdot x^p,$ and pass from $\k$ to $\C,$ by Lefschetz principle, \S\ref{Sec.Prelim.Lefschetz}.
\bee[Step 1.]
\item
 The polynomial $a(x)$ has at least two distinct roots. Split the fraction $\frac{1}{a(x)}$ into simple fractions,
  $\sum\frac{c_{ij}}{(x-z_j)^i}.$ We prove: $\frac{1}{a(x)}$  has at least two non-zero residues, i.e. in this decomposition there are at
   least two linear terms, $\frac{1}{x-z_j}$.

\item The power series $y(x)$ must satisfy the equation of Lemma \ref{Thm.Exp.Der.Case.n=1}.
  We prove:  $y(x)$ must also satisfy the equation $\sum \la_j \cdot ln\frac{z_j-y}{z_j-x}=const.$

\item We prove: the algebraic function $y(x)$ cannot satisfy this transcendental equation.
\eee

\medskip

\bee[\!\!\!\bf Step 1.]

\item
  The polynomial $a(x)$ has at least two distinct roots. Factorize it, $a(x)=c\cdot \prod^r_{j=1}(x-z_j)^{m_j},$ with $r\ge2$ and $z_i\neq z_j.$
   Accordingly split the fraction,
  \beq
\frac{1}{a(x)}=\sum^r_{j=1} \sum^{m_j}_{i=2}\frac{b_{j,i}}{(x-z_j)^i}+\sum^r_{j=1}\frac{\la_j}{x-z_j},\quad
 \text{for some } \{b_{j,i}\}, \{\la_j\} \text{ in }\C.
  \eeq
 We claim: at least two residues of $\frac{1}{a(x)}$ are non-zero, i.e. at least two of $\la_j$'s do not vanish.
 First we compute the total residue, $\sum_j res_{z_j}\frac{1}{a(x)}=res_\infty \frac{1}{a(x)}=0.$ (Because $deg\ a(x)\ge2.$)
  Thus it is enough to prove: at least one of $\la_j$ is non-zero.

\medskip

  Suppose all the residues vanish, then the antiderivative of $\frac{1}{a(x)}$ is a rational function, presentable as $\frac{p(x)}{\prod (x-z_j)^{m_j-1}}.$
   When $x\to \infty$ this function   decays with the rate $-deg(p(x))+\sum_j (m_j-1)\le (\sum_j  m_j)-r.$
    But then the decay rate of $\frac{1}{a(x)}$  at $x=\infty$
    is at most $deg\ a(x)-r+1< deg\ a(x) .$ Hence the contradiction.

\item
 Denote $y(x)\!=\!e^\xi x,$ then
 $y(x)\!\in\!\k\bl x\br,$ by Corollary \ref{Thm.If.D^infty.then.algebraic}.
 However, by Lemma \ref{Thm.Exp.Der.Case.n=1}, $y(x)$ satisfies the equation
 \beq\label{Eq.inside.proof.Algebraicity}
Q(x)-Q(y)+\sum^r_{j=1}  \la_j\cdot ln \frac {y-z_j}{x-z_j} =1,\qquad \text{ for some }0\neq Q(x)\in \C(x).
 \eeq
We claim: in this case $y(x)$ satisfies two equations:
\beq
Q(x)-Q(y)=const \quad  \text{ and }\quad  \sum^r_{j=1}  \la_j\cdot ln \frac {y-z_j}{x-z_j} =const.
\eeq
\bei
\item {\em Case 1,} the coefficients $\{\la_j\}$ are pairwise $\Q$-dependent, i.e. $\la_j\in \Q\la_i$ for all $i,j.$
 Then, by scaling equation
 \eqref{Eq.inside.proof.Algebraicity}, we can assume: $\{\la_j\in \Z\}.$ The equation becomes $Q(x)-Q(y)+ln(\prod_j q_j(x,y)^{\la_j})=const,$
  where $q_j(x,y)^{\la_j}\in \C(x,y).$

  And then the claim follows, because $ln$ is a transcendental function.

\item {\em Case 2,} the coefficients $\{\la_j\}$ are pairwise $\Q$-independent, i.e. $\la_j\not\in \Q\la_i$ for all $i\neq j.$
 Pass to the differential forms, $d(Q(x)-Q(y(x)))+\sum_j \la_j \frac{dq_j(x,y(x))}{q_j(x,y(x))}=0.$ Here $d(Q(x)-Q(y(x)))$ is an exact differential
  form on the compact Riemann
  surface of the algebraic function $y(x).$ Therefore its residues at all the points of this Riemann surface (denote it $X$) vanish.
   Therefore
   \beq
   \sum_j\la_j \cdot res_p \frac{dq_j(x,y(x))}{q_j(x,y(x))}=0\qquad \forall p\in X.
   \eeq
Finally observe: $res_p \frac{dq_j(x,y(x))}{q_j(x,y(x))}=ord_p q_j(x,y(x))\in \Z.$ Therefore for each $p\in X$ we get a $\Z$-linear
  combination of $ \la_\bullet$
 that sums up to zero. As all $ \la_\bullet$ are $\Q$-linearly independent, we get: $ord_p q_j(x,y(x))=0$ $ \forall p\in X.$ Thus each $q_j(x,y(x))$
  is a meromorphic function without zeros and poles on the compact Riemann surface $X.$ Therefore $q_j(x,y(x))=const$ for each $j.$

\item {\em Case 3.} Consider the $\Q$-vector space, $Span_\Q(\la_1,\dots ,\la_r).$ After a permutation of indices we can assume: its
 basis is $\la_1,\dots,\la_l,$ for some $l\le r.$ Suppose $l<r,$ then $\la_{l+1},\dots,\la_r\in Span_\Q(\la_1,\dots ,\la_l).$
  Therefore $ \sum^r_{j=1}\la_j \cdot ln(q_j(x,y))= \sum^l_{j=1}\la_j(\sum_i a_{i,j} \cdot ln(q_j(x,y)),$ for some $a_{i,j}\in \Q.$
   Scaling the equation we can assume  $a_{i,j}\in \Z.$ Then the equation becomes $Q(x)-Q(y)+\sum^l_{j=1}\la_j\cdot ln(\prod_i q_j(x,y)^{a_{ij}})=const.$
    And now invoke Case 2.
  \eei

\item The power series $y(x)=e^\xi x=x+x^p\cdot(\dots),$  with $p\ge2,$ satisfies the equation $F(x,y)=\sum_j \la_j\cdot ln\frac{z_j-y}{z_j-x}=const.$
 We expand $F(x,y(x))$ into the power series:
 \bei
 \item For the roots $z_j\neq0$ we have:  $ln\frac{z_j-y(x)}{z_j-x}=ln(1+\frac{x-y(x)}{z_j-x})\in (x)^p.$
 \item For the root $z_j=0$ (if this term exists) we have:  $ln\frac{ y(x)}{ x}= x^{p-1}\cdot (\dots).$
 \eei
Therefore the equation is: $F(x,y)=0.$ And moreover, if $z_j=0$ then $\la_j=0.$ (By comparison of the $x$-orders.)

Therefore all the constants $z_\bullet$ in $F(x,y)$ are non-zero. Expand $F(x,y(x))$:
\beq
0=\sum_j \la_j\cdot ln\frac{z_j-y(x)}{z_j-x}=\sum_j \la_j\cdot ln(1-\frac{y(x)-x}{z_j-x})=\hspace{2cm}
\eeq
$$
\hspace{5cm}=(y(x)-x)\cdot \Big[-L(x)+L'(x)\cdot\frac{y(x)-x}{2}-
 L''(x)\cdot\frac{(y(x)-x)^2}{3}+\cdots\Big].
$$

Here $L(x):=\sum_j\frac{\la_j}{z_j-x}\in \k[\![x]\!],$ thus $L\neq0$ and $\ord_x L(x)<\infty.$ But then
\beq
\ord_x L(x)<\min\{\ord_x L'(x)\cdot\frac{y(x)-x}{2},\ L''(x)\cdot\frac{(y(x)-x)^2}{3},\dots\}.
\eeq
Hence the contradiction.
\epr
\eee

\subsection{An extension to the case $n\ge2.$ Computation of $e^\xi x$ for derivations with separated variables}\label{Sec.Exp.Der.computation.Separated.Derivations}
As before we transform the (transcendental) computation of $e^\xi x$ into the solution of implicit function equation.

Let $\k$ be a local domain and $x\!=\!(x_1,\dots,x_n).$  Fix a derivation $\xi\!=\!\sum \xi_i \di_{x_i}\in \Der^{nilp}_\k(\k[\![x]\!]).$
 We can assume (wlog): $\xi_i\neq0$ for all $i=1\dots n.$ Indeed, if $\xi_n=0$ then $\xi$ is $x_n$-linear.
  Then we  pass to the $\K$-algebra $\K[\![x_1\dots x_{n-1}]\!],$
  where $\K=\k[\![x_n]\!].$

Suppose $\xi$ satisfies: $\frac{\xi_i}{\xi_j}=\frac{\ta_j(x_j)}{\ta_i(x_i)}$ for all $i,j,$ and some $\ta_i\in \Frac(\k[\![x_i]\!]).$
   Then 
\beq\label{Eq.Derivation.Separated}
\xi=b\cdot \sum^n_{i=1}\frac{1}{ \ta_i(x_i)}\frac{\di}{\di x_i},  \quad  \text{ for some  }
  \quad      b\in \Frac(\k[\![x]\!]).
\eeq

This is a rather special class of derivations. Yet they show up in numerous applications.

\medskip

The element $ \ta_i(x_i)$ is a   Laurent power series  with finite negative part.
 Take its residue, $\la_i:=res_o\ta_i(x_i)\in  \Frac(\k),$ i.e. the coefficient of the monomial $\frac{1}{x_i}.$

Take the  ``adjusted order of the pole", $p_i+1:=-ord_{x_i}[\ta_i(x_i)-\frac{\la_i}{x_i}].$ Thus $p_i\in \Z\smin \{0\}.$

\bel \label{Thm.Exp.Der.Separated.Variables}
\bee
\item
The derivation $\xi$ of \eqref{Eq.Derivation.Separated} can be presented in the form
\beq\label{Eq.IF.eq.def.exp.raw}
\xi=b \cdot \sum_{i=1}^n\frac{1}{c_i\cdot \di_{x_i} \frac{1}{ a_i(x_i)^{p_i} }+\frac{\la_i}{x_i}}\frac{\di}{\di x_i},
\eeq
where:  \ $ b \in \Frac(\k[\![x]\!]),$   \   $c_i,       \la_i \in \Frac(\k),$ and $a_i(x_i)-x_i \in (x_i)^2 \sset \Frac(\k)[\![x_i]\!].$

      The solution of the ODE $x'=\xi,$ $x(o)=x^o$  satisfies:
             \beq\label{Eq.IF.eq.defining.integral.curve}
              \frac{ c_1}{ a_1(x_1)^{p_1}}- \frac{ c_1}{ a_1(x^o_1)^{p_1}}+\la_1\cdot ln\frac{x_1}{x^o_1}
 =\cdots=             \frac{ c_n}{ a_n(x_n)^{p_n}}- \frac{ c_r}{ a_n(x^o_n)^{p_n}}+\la_n\cdot ln\frac{x_n}{x^o_n}.
\eeq

         \item Denote $y_i =e^\xi x_i.$ Then   $y_1\dots y_n$ are the unique solutions of the  implicit function equations
       \beq\label{Eq.IF.eq.defining.exp}
       1=\int^{y_i}_{x_i}\frac{c_i\cdot \di_i \frac{1}{ a_i(s_i)^{p_i}} +\frac{\la_i}{s_i}}{b(x(s_i,x^o))}ds_i,\quad \quad i=1\dots n  .
       \eeq
      (Here the elements  $x(s_i,x^o) $ are determined by \eqref{Eq.IF.eq.defining.integral.curve}.)
\\
In particular,
\bei
\item  If $b\in \Frac(\k),$ then   the equations are
 $   \frac{c_i}{ a_i(y_i)^{p_i}} -   \frac{c_i}{ a_i(x_i)^{p_i}}  +\la_i\cdot ln\frac{y_i}{x_i}=b$ for $i=1\dots n.$
 \item  Moreover, if $\la_i\!=\!0,$  then
 $e^\xi x_i\!=\!a^{-1}_i[\frac{a_i(x_i)}{\big(1+\frac{b}{c_i} \cdot a_i(x_i)^{p_i}\big)^\frac{1}{p_i}}].$  Here $a^{-1}_i(\dots)$ is the inverse of $a_i(x_i)$.
\eei \eee
\eel
\bpr The question is essentially reduced to the case $n\!=\!1.$
\bee
\item
Expand $\ta_i(x_i)$ of \eqref{Eq.Derivation.Separated} into Laurent power series in $x_i.$
  It has a finite singular part. Denoting $\K=\Frac(\k)$ we can present:
\beq
\int (\ta_i(x_i)-\frac{\la_i}{x_i})dx_i=\frac{1}{c_i}(x_i^{-p_i}+\cdots)\in Span_\K(\frac{1}{x^{p_i}_i},\dots,\frac{1}{x_i})+
 \K[\![x_i]\!].
\eeq
 Here $c_i\in\K  .$
Rewrite this as:
  $\int (\ta_i(x_i)-\frac{\la}{x_i})dx_i=\frac{1}{c_i\cdot a_i(x_i)^{p_i}}.$
 This gives equation \eqref{Eq.IF.eq.def.exp.raw}.

Consider now the Pfaffian system $\frac{dx_1}{\xi_1}=\cdots =\frac{dx_r}{\xi_r},$ where
  $\xi_i=\frac{b}{c_i\cdot \di_{x_i} \frac{1}{ a_i(x_i)^{p_i} }+\frac{\la_i}{x_i}}.$ It is integrated directly, giving equation \eqref{Eq.IF.eq.defining.integral.curve}.

\item Substitute $x_i=x^o_i\cdot (1+\De_i)$ into \eqref{Eq.IF.eq.defining.integral.curve}. Note that for all $i,j$ the corresponding equations contain term linear
 in $\De_i,\De_j.$ Thus we have the solution $\De_i=f(\De_j,x^o_i,x^o_j)\in \Frac(\k[\![x^o_i,x^o_j]\!])[\![\De_j]\!].$
  Substitute this into \eqref{Eq.IF.eq.def.exp.raw}. Then $\frac{dx_i}{dt}=\xi_i$ is transformed into $dt=\frac{dx_i}{\xi_i},$ where $\xi_i$ depends on $x_i,x^o$ only. Finally integrate this.
\epr
\eee

\bex
\bee[\bf i.]
 \item The particular case   $\xi=\xi_1(x_1,\dots,x_n)\cdot  \frac{\di}{\di x_1} $
  is   one-dimensional, over the $\K$-algebra   $\K[\![x_1]\!],$ where  $\K=\k[\![x_2,\dots,x_n]\!].$

 E.g. consider $\xi=b\cdot x_1\frac{\di}{\di x_1},$ with $b\in (x_2,\dots, x_n)\sset \k[\![x_2,\dots,x_n]\!].$
 Here $e^\xi x_1=e^{b }x_1.$
\item
 Another essentially one-dimensional case is
  $\xi=\sum \xi_i(x_i)\frac{\di}{\di x_i},$ where $\xi_i(x_i)\in (x_i)^2\sset   \k[\![x_i]\!] .$
   Then $e^{\xi}x_i=e^{\xi_i(x_i)\frac{\di}{\di x_i}}x_i.$  The results/examples of \S\ref{Sec.Exp.Der.n=1} readily apply.

\item
Consider $\xi=x^{p_1}_1x^{p_2}_2\cdot (\frac{1}{p_1 x^{p_1-1}_1}\di_{x_1}+\frac{1}{p_2 x^{p_2-1}_2}\di_{x_2}),$ where $p_1,p_2\ge 1.$
 Here equation \eqref{Eq.IF.eq.defining.integral.curve}  is: $x^{p_1}_1-(x^o_1)^{p_1}= x^{p_2}_2- (x^o_2)^{p_2}.$ Then
 equation \eqref{Eq.IF.eq.defining.exp} gives for $i=1:$
\beq
  1=\int^{y_1}_{x_1}   \frac{\di_{s_1} s^{p_1}_1\cdot ds_1}{s^{p_1}_1(s^{p_1}_1-x^{p_1}_1+x^{p_2}_2) }  \stackrel{t=s^{p_1}_1}{=}
  \frac{1}{ x^{p_1}_1-x^{p_2}_2}\cdot ln\big[1-\frac{ x^{p_1}_1-x^{p_2}_2}{t}\big]\Big|^{t=y^{p_1}_1}_{t=x^{p_1}_1}.
\eeq
We get the equation for $y^{p_1}_1,$ which is easily solvable:
\beq
e^\xi x_1=y_1(x)=\frac{x_1}{\sqrt[p_1]{1-x^{p_2}_2\cdot \frac{ e^{x^{p_1}_1- x^{p_2}_2}-1}{x^{p_1}_1-x^{p_2}_2}}}\in D^2(\k[x_1,x_2]).
\eeq
The formula for $e^\xi x_2$ is similar.
\eee
\eex

\section{$e^\xi$ acts on various rings of power series over normed fields}
Fix  a field $\k\supseteq \Q$   with a (non-trivial) norm  $|\dots|$.  (The classical examples are $,\Q,\R,\C$ and p-adics)
 Let $x=(x_1,\dots,x_n)$ and denote $|\um|=\sum _i m_i.$
 To compare (formal) power series, we write $f\preceq \tf,$ for $f(\ux)=\sum c_\um \ux^\um$ and $\tf(\ux)=\sum \tc_\um \ux^\um,$
  if $ \sum_{|\um|=d}|c_\um|\le  \sum_{|\um|=d}|\tc_\um|$ \  for each $1\ll d\in\N .$

Take a  (not necessarily Noetherian) differential $\k$-algebra $\k[\ux]  \sset R\sset \k[\![\ux]\!]$ that satisfies the following conditions.
\bee[\bf a.]
\item ($R$ is closed under majorization) If $\tf\in R$ and $f\preceq   \tf,$   then $f\in R;$
\item ($R$ admits the implicit function theorem in two variables) If $f(x,y)\in R$ with $f(0,0)=0$ and $\di_y f|_{(0,0)}\neq0,$ then the unique
 formal solution $f(x,y(x))=0,$ $y(x)\in (x)\sset \k[\![x]\!],$ belongs to $R.$
\eee
The main examples are rings with coefficients of controlled growth,
\beq
 R^{\{\la_\bullet\}}=\{\sum c_\um \ux^\um| \quad \sum_{|\um|=d}|c_\um|\le\la_d,\ \forall d\}\sset\k[\![x]\!].
\eeq
Here  $\{\la_\bullet\}$ are  some prescribed sequences, depending on $\sum c_\um \ux^\um.$     See Example \ref{Ex.Exp.Der.acts.for.normed.field}.

 \bthe\label{Thm.Exp.Der.Bounded.Coefficients}
   $e^\xi $ acts on $R $ for any $\xi\in   \Der^{nilp}_\k(R).$
 \ethe
\bpr
 {\bf Preparations.}
 \bei
 \item We claim: $R$ is local, i.e. for any $h(x)\in R$ with $h(o)=0,$ we have $\frac{1}{1+h(x)}\in R.$
  Indeed, consider the equation $f(x,y)\!:=(1+h(x))\cdot  y-1\!=\!0.$ Here $\di_y f(x,y)|_{o,o}\neq0$  and the formal solution is
    $y\!=\!\frac{1}{1+h(x)}\in \k[\![x]\!].$
     By the assumption
   $b.$ we get: $\frac{1}{1+h(x)}\!\in\! R.$

 \item  We claim: the ring $R$ is closed under integration,  i.e. if $f(\ux)\in R$ then $\int f(\ux)dx_i\in R$ for each $i.$
 Indeed, let $f(x)=\sum c_\um \ux^\um,$ then  $\int f(\ux)dx_i\preceq x_i\cdot f(\ux)\in R.$

In particular, $ln(1+x)=\int\frac{dx}{1+x}\in R.$

\item (The diagonal reduction) For every $f(x)\in \k[\![x]\!]$ take $f(x_1,\dots,x_1)\in \k[\![x_1]\!].$
 Observe: $f(x_1,\dots,x_1)\preceq f(x).$ Therefore $f(x_1,\dots,x_1)\in R.$

\item (The $\Q_{\ge0}$-majorization of $f$) Take any $f(x)=\sum_\um c_\um \ux^\um\in \k[\![x]\!].$ Fix some non-negative rationals $c_{\Q,\um}\in \Q_{\ge0}$ satisfying:
 $|c_\um|\le |c_{\Q,\um}|\le 2|c_\um|.$ Denote $f_{\Q_{\ge0}}:=\sum_\um c_{\Q,\um} \ux^\um.$ Observe: $f\preceq f_{\Q_{\ge0}} \preceq 2f.$
  Therefore $f\in R$ \iff $f_{\Q_{\ge0}}\in R.$

\item ($R$ admits compositions, \S\ref{Sec.Prelim.Notations.Conventions}) Let $f(x),g(x)\in R,$ with $g(0)=0.$ It is enough to prove $f(g(x))\in R$
 for the case: the coefficients of $f(x),g(x)$ all belong to $\Q_{>0}.$ Wlog we can assume $f'|_{x=0}\neq1.$

 Consider the implicit function equation $y=g(x)+f(y).$ Its unique solution is $y(x)=g(x)+f(g(x)+f(\cdots))\in R.$ Thus $f(g(x))\preceq y(x).$
 Hence $f(g(x))\in R.$

\item (Eliminating the linear part of $\xi$) Expand the derivation, $\xi=\xi_1+\xi_{\ge2}\in  \Der^{nilp}_\k(R),$ where $\xi_{\ge2}\in (x)^2\cdot \Der_\k(R),$
 while $\xi_1=\sum^n_{i=1} l_i(x)\di_i,$ with $l_\bullet$ linear forms. As $\xi$ is $(x)$-adically nilpotent, $\xi_1$ is $(x)$-adically nilpotent as well.
  But then $\xi_1$ is just locally nilpotent on $\k[x].$  Thus the operator $e^{\xi_1}$ acts on $\k[\![x]\!]$ by a $GL(n,\k)$-transformation.
   Therefore, to prove that $e^\xi x\in R,$ it is enough to prove that $e^{-\xi_1}\cdot e^\xi x\in e^{-\xi_1} R.$
 Note that the subring $e^{-\xi_1} R\sset\k[\![x]\!]$ satisfies the assumptions a., b.
    Finally,
    by the Baker-Campbell-Hausdorff formula, we have $e^{-\xi_1}\cdot e^\xi=e^\txi,$ where $\txi\in (x)^2\cdot \Der_\k(R). $

    Therefore below we assume:  $\xi\in (x)^2\cdot \Der_\k(R). $
\eei

\medskip {\bf Reduction to the case $\pmb{n=1.}$}
We want to verify: $e^\xi x\in R.$

Let $\xi=\sum_i \xi_i(x)\frac{\di}{\di x_i}.$ To each $\xi_i$ associate its $\Q$-majorization, $\xi_{\Q_{\ge0},i}\in R,$ as above.
 Thus $\xi_i\preceq \xi_{\Q_{\ge0},i}\le 2\xi_i.$ Take the diagonal reduction (as above)
  $\txi  := (\sum_i   \xi_{\Q_{\ge0},i}(x_1,\dots,x_1)\big)\di_{x_1} .$ Thus $\txi\in (x_1)^2\cdot \Der_\k(R).$

 Observe:  $ e^\xi x_i \preceq e^{\txi }x_1$ for each $i.$ Therefore it is enough to prove the statement for the subring  $  R|_{x_1=\cdots=x_n}\sset R,$
 i.e. for the one-dimensional case.

\medskip
   {\bf The case of $\pmb{n=1.}$}
  Denote $y(x)=e^\xi x$ and invoke Lemma \ref{Thm.Exp.Der.Case.n=1}. Thus $y(x)$ satisfies:
\beq
\frac{1 }{c\cdot y^p\cdot (1+\ta(y)) }-\frac{1 }{c\cdot x^p\cdot (1+\ta(x)) }+\la\cdot \ln(\frac{y}{x})=1,
\eeq
 where $ p\ge1,$ $ \ta(o)=0,$ and $y^p\cdot (1+\ta(y))\in \int R=R.$

    Expand $y=x\cdot (1+\ty).$   Then this  equation becomes
\beq
\frac{1 }{c\cdot (1+\ty)^p\cdot (x\cdot (1+\ty)) }-\frac{1 }{c\cdot    (1+\ta(x)) }+\la\cdot \ln(1+\ty)=1.
\eeq
Clearing the denominators we get an implicit function equation, $R_{x,\ty}\ni F(x,\ty)=0.$ It satisfies   $\di_\ty F|_{o,o}\neq0,$ and its (unique)
  formal solution is $\ty(x)=\frac{e^\xi x-x}{x}\in \k[\![x]\!].$
   Thus $\ty(x)\in R.$ Therefore $e^\xi x\in R.$
\epr

\bex\label{Ex.Exp.Der.acts.for.normed.field}
    Fix an eventually  log-convex sequence $\la_\bullet\sset \R_{\ge1}$, i.e. the sequence $\{\frac{ \la_{d+1}}{ \la_d}\}$ is eventually increasing.
 Consider the (formal) Denjoy-Carleman class (in $n$ variables)
 \beq
 R^{\{\la_\bullet\}}:=\{f(\ux)=\sum c_\um \ux^\um\big|\ \exists C_f,\rho_f\ \text{ such that }
   |c_\um|\le C_f\cdot (\rho_f)^{|\um|}\cdot \la_{|\um|}, \ \forall \um \}\sset \k[\![\ux]\!],
 \eeq
 see e.g. \cite{Thilliez2008}. To apply lemma \ref{Thm.Exp.Der.Bounded.Coefficients} to this ring we verify the assumptions.
\bee[\quad \bf a.]

\item   Suppose $f\preceq \tf\in R^{\{\la_\bullet\}}.$ Thus $\sum_{|\um|=d}|c_\um|\le\sum_{|\um|=d}|\tc_\um|.$
 Therefore (for each $\um$) $|c_\um|\le \bin{n+|\um|}{|\um|}\cdot c_\tf\cdot \rho_\tf^{|\um|}\cdot \la_{|\um|}.$
  It remains to observe:
  \beq
  \bin{n+d}{d}\le (1+d)\cdots (n+d)\stackrel{d\ge n}{<}(2d)^n< M^d \text{ for }M\gg1 .
\eeq
Thus $f\in R^{\{\la_\bullet\}}.$

\item   The ring $R^{\{\la_\bullet\}}$  satisfies the Implicit Function Theorem \cite{Komatsu}.

\eee
\noindent   Therefore  $e^\xi $ acts on $R^{\{\la_\bullet\}}$ for any $\xi\in \Der^{nilp}_\k(R^{\{\la_\bullet\}})$.

\medskip
The simplest particular cases of such rings (where all $x$-adically nilpotent derivations exponentiate) are:
\bei
\item For $\la_d=1$ one has  $R^{\{\la_\bullet\}}=\k\{\ux\}$. Thus $e^\xi $ acts on the ring of $\k$-analytic power series.
 This is classic, both for $\R,\C$ and for p-adics.

\item Take $\la_d= (d!)^\al$ for $\al>0.$ (This is called  a Gevrey sequence.)
 Then $R^{\{\la_\bullet\}})$ is related to Gevrey $G^{1+\al}$ regularity class,
    \cite[Example 2.1.2]{Thilliez2008}.

In this case the action of $e^\xi $ on Gevrey power series is also known, e.g. \cite{Cano}.

\item In numerous other cases of log-convex sequences we get the action    $e^\xi\circlearrowright R^{\{\la_\bullet\}}.$
\eei

 \eex

\section{$e^\xi $ does not act on the ring $\quots{C^\infty(\R^n,o)}{J}$}\label{Sec.Cinfty.case}
An ideal $J\sset (x)\sset C^\infty(\R^n,o)$ defines the closed subgerm $o\in V(J)\sseteq (\R^n,o).$
 Below we assume: $V(J)\neq \{o\},$ i.e. any   representative of $V(J)$ contains more than one point.

 Take a derivation $\xi\in \Der_\R(\quots{C^\infty(\R^n,o)}{J}).$
 Present it by   $\xi=\sum^n_{i=1} \xi_i(\ux)\frac{\di}{\di x_i}\in \Der_\R(C^\infty(\R^n,o))$ such that $\xi(J)\sseteq J.$
 We assume $\xi|_{V(J)}\neq0$, i.e. at least one of $g_i(\ux)$ does not vanish identically on the set-germ $V(J)$.

\bthe\label{Thm.C.infty.case.exp.does.not.act}
 If $\xi|_{V(J)}\neq 0$
 then  the operator $e^\xi $  does not act on $\quots{C^\infty(\R^n,o)}{J}$.
\ethe
Namely, for a given $\xi$ we construct a function-germ $f\in C^\infty(\R^n,o)$ such that $e^\xi f$ diverges on a convergent sequence of points $o\leftarrow \{z_\bullet\}\sset V(J)$.
\bpr
It goes in several steps. In steps 1-3 we work in the ring $C^\infty(\R^n,o)$.
\bee[\em Step 1.]
\item (Preparations) If $e^\xi$ acts on  $C^\infty(\R^n,o),$ then necessarily $\xi(x)\sseteq (x),$ i.e. $\xi|_o=0.$

Furthermore,    the case of $C^\infty(\R^n,o)$ is reduced to the case $C^\infty(\R^1,o).$
\item If  $\xi$ has one-sided-isolated zero at $o\in \R^1,$ then  $e^\xi$ does not act on $C^\infty(\R^1,o).$
\item  $e^\xi$ does not act on $C^\infty(\R^1,o) $   for an arbitrary $\xi$.
\item Extension   to the ring $\quots{C^\infty(\R^n,o)}{J}.$
\eee
\bee[\!\!\!\!\bf\text{Step} 1.\!]
\item (Preparations)
\bee[\hspace{-0.6cm}{$\bullet$}]
\item Suppose $\xi|_o\neq0,$ i.e. $(\xi_1(o),\dots,\xi_n(o))\neq (0,\dots,0).$
 Such a vector field can be rectified, \cite[\S7]{Arnol'd}. Namely, there exists a $C^\infty$-coordinate change, $\Phi\circlearrowright (\R^n,o),$ $x\to \tx(x),$
  satisfying: $\Phi_* \xi=\frac{\di}{\di x_1}.$ (See \S\ref{Sec.Prelim.Derivations}.iii.)   Thus we can assume
  $\xi=\frac{\di}{\di x_1}.$

  Then $e^\xi$
 does not act on $C^\infty(\R^n,o).$ E.g. take any analytic series $f(x_1)=\sum c_i x^i,$ with $\sum c_i=\infty.$
  Then  $e^\xi f|_o =\infty.$

Therefore below we assume $\xi|_o=0,$ i.e. $\xi(x)\sseteq (x).$

\item (Reduction to one-dimensional case) Take the $\hx_1$-axis, $\phi: \R \hx_1 \hookrightarrow \R^n.$
Let $\xi=\sum \xi_i \frac{\di}{\di x_i}\neq0,$ we can assume (for generic coordinates) $\phi^*(\xi_1)\neq0,$ thus the
 restricted component of the vector field,
 $\xi|_{\R\hx_1}:= \phi^*(\xi_1) \frac{\di}{\di x_1},$ is non-zero.

Fix a sequence of points on $\hx_1$-axis,  $(z_\bullet,0\dots,0)\in \R^n,$ such that $z_\bullet\to o\in \R^1.$ Below we construct
 $f(x_1)\in C^\infty(\R^1,o)$ such that $e^{\xi|_{\R\hx_1}}f$ diverges at the points $z_\bullet.$
 Consider $f(x_1)$ as a function on $(\R^n,o)$ via the projection $(\R^n,o)\to(\R^1_{x_1},o).$
 We want to deduce:
  $e^\xi f(x)$ diverges at these points as well.  It is enough to assume:
\beq\label{Eq.Cinfty.case.reduction.to.dim=1}
\forall\ j\ge d\ge1 :\quad\quad\quad
\big|\xi^j
f |_{z_d}- (\xi|_{\R\hx_1} )^j f |_{z_d}\big| <
\frac{1}{2}\big|(\xi|_{\R\hx_1} )^j f |_{z_d}\big|.
\eeq
Indeed, in  this case $ e^{\xi|_{\R\hx_1} } f |_{z_d}  =\infty$ will imply  $\big|e^\xi f|_{z_d}\big|\ge \frac{1}{2}\big|e^{\xi|_{\R\hx_1} } f |_{z_d}\big|=\infty.$

\quad Observe, the term $f  - (\xi|_{\R\hx_1} )^j f $ contains the derivatives $\di^\bullet_1 f$ of order$\le j-1$ only.
 Therefore, to ensure \eqref{Eq.Cinfty.case.reduction.to.dim=1}, it is enough to show: for the chosen sequence $\{z_\bullet\}$ and {\em any}
  fast growing sequence of numbers $\{N_\bullet\}$ there exists $f\in C^\infty(\R^1,o)$ satisfying:
\beq\label{Eq.fast.growing.derivatives}
 \forall\ j\ge d\ge1:\quad\quad\quad
 (\xi|_{\R\hx_1})^{j+1}f|_{z_d}> N_j\cdot (\xi|_{\R\hx_1})^j f|_{z_d}.
\eeq
 (This condition ensures, in particular: $e^{\xi|_{\R\hx_1}}f|_{z_d}=\infty$.)

 Altogether, the question is reduced to the ring $C^\infty(\R^1,o),$ and it remains to construct such a function $f.$
\eee

\medskip

\item Below $\xi=g(x)\frac{d}{dx}.$
Consider first the special case: $g(x)\in C^\infty(\R^1,o)$ has a
one-sided-isolated zero at $o$.
 Namely, $g(0)=0$ and for some (small) $\ep>0$ holds:  $g$ is smooth on $(o,\ep)$   and $g(x)>0$ for $x\in (0,\ep]$.
\bee[\hspace{-0.8cm}\bf\text{Part} A.]
\item
 We want to eliminate the coefficient $g(x)$ by changing the variables, $g(x)\frac{d}{dx}=\frac{d}{dy}.$ Thus
  consider the ODE $\frac{dy}{dx}=-\frac{1}{g(x)}$ on $(o,\ep)$.
 Its   solution is   $y(x)=\int_x^{\ep} \frac{d\tx}{g(\tx)}\in C^\infty(o,\ep).$
  It satisfies:
 \beq
 \liml_{x\to 0^+} y(x)=\infty,\quad \quad
 y >0\ on\ (o,\ep),\quad\quad
   (o,\ep)\isom{y(x)}(y(\ep),\infty) \text{ is a diffeomorphism.}
 \eeq

 Therefore, for $x\in (o,\ep)$ we present  $e^{g(x)\frac{d}{dx}}=e^{- \frac{d}{dy}}$ and it is enough to construct a function $\tf(y)\in C^\infty(y(\ep),\infty)$ satisfying:
\bee[\bf i.\!]
\item
 $\tf(y)$ is flat at $\infty$, i.e. for each $j\in \N$ holds: $\liml_{y\to \infty}\frac{d^j \tf(y)}{dy^j}=0$.
\item
There exists a sequence of points $y_\bullet\to \infty$ such that the series
 $e^{d/dy}\tf|_{y_d}=\sum^\infty_{j=0}\frac{1}{j!}\frac{d^j\tf}{dy^j}|_{y_d}$ diverges for each $d$.
\item For the   sequence of numbers $\{N_\bullet\}$ of equation \eqref{Eq.fast.growing.derivatives} the bound
  holds:
  \beq
   \frac{d^{j+1}\tf}{dy^{j+1}}|_{y_d} > N_j\cdot \frac{d^j\tf}{dy^j}|_{y_d}\quad \text{ for }j\ge d.
  \eeq
\eee

Once we have such a function, $\tf(y)$, we define
\beq\label{Eq.inside.proof.smooth.case}
f(x):=\Big\{\ber
\tf(y(x)),\ x\in (o,\ep)\\0,\   x\le0\eer.
\eeq

 It is $C^\infty(-\ep,\ep)$ by construction and flat at $x=0$.
 In addition
 \beq
 (e^{g(x)\frac{d}{dx}}f(x))|_{z_d}=(e^{-\frac{d}{dy}}\tf(y))|_{y(z_d)}=\infty, \quad \forall d\in\N.
 \eeq
Thus
 $e^{g(x)\frac{d}{dx}} f(x) \not\in C^\infty(\R^1,o)$. Finally, $\big(g(x)\frac{d}{dx}\big)^j f|_{z_d}$ grows strongly,
  so that condition \eqref{Eq.fast.growing.derivatives} is satisfied.

\

The case $g(x)<0$ on $(o,\ep)$ is done similarly, with $\frac{dy}{dx}=\frac{1}{g(x)}.$

\item  We construct $\tf(y)\in C^\infty(\R_{\ge0})$   in the form $\tf(y)=\frac{\tau(\{y\})}{\lfloor y\rfloor+1}.$
 Here $\{y\}$ is the fractional part, $\lfloor y\rfloor$ is the integer part.
 The function $\tau\in C^\infty([0,1])$ is assumed flat at the points $0$ and  $1$.
 Thus:
 \bei
 \item  $\tf(y)$ is flat on the set $\N$ and $\tf\in C^\infty(\R_{\ge0})$.
\item    $\liml_{y\to \infty}f^{(j)}(y)=0$, because   $\lfloor y\rfloor$ is constant on $\R_{\ge0}\smin\N$ for any fixed $j$.
 Thus $\tf(y)$ is flat at $ \infty$.
\eei

Now we apply the exponential:
\beq
e^{-\frac{d}{dy}}\tf(y)=\sum^\infty_{j=0}\frac{(-1)^j}{j!}\frac{d^j\tf(y)}{dy^j}=\Bigg\{
\ber\ 0,\ y\in \N\\
\frac{1}{\lfloor y \rfloor+1}\sum^\infty_{j=0}\frac{(-1)^j}{j!}\frac{d^j\tau(\{y\})}{dy^j}\eer,
\ y \in\R_{\ge0}\smin \N.
\eeq
Finally we assume that the flat
function $\tau(y)$ has fast growing derivatives  for at least one
point $y_0\in(0,1)$:
\beq\label{Eq.inside.proof.Step.3.part.B}
\frac{d^{j+1}\tau(y)}{dy^{j+1}}|_{y_0} > j!\cdot p_j\cdot
\frac{d^j\tau(y)}{dy^j}|_{y_0}\gg1. \eeq
 Then
 $e^{\frac{d }{dy }}\tf(y)$ diverges at each point of the sequence $ \{y_0\}+\N$, for $y_0\not\in\Z$.
\eee
\medskip

Now, for $f$ as in equation \eqref{Eq.inside.proof.Step.3.part.B},
  $e^{g(x)\frac{d}{dx}}f(x)$ diverges on  the sequence of
points $z_\bullet:= x(\{y_0\}+\N)\to0$.
 Thus $e^{g(x)\frac{d}{dx}}f(x)\not\in C^\infty(\R^1,o)$.

  Moreover,
    the sequence $\big(g(x)\frac{d}{dx}\big)^j f|_{z_\bullet}$ is fast growing in the sense of equation \eqref{Eq.fast.growing.derivatives}.
 Thus also for the initial  derivation, $\xi=\sum a_i \frac{\di}{\di x_i}\in \Der_\R(C^\infty(\R^n,o))$,
  we have the divergence, $e^{\xi}f(x)\not\in C^\infty(\R^n,o)$.

\

\item
Finally consider the general case: $g\neq 0$, but the zero of $g(x)$ at $x=0\in \R^1$ is not necessarily one-sided-isolated.
  Then there exists a sequence of points
  $\R^1\supset z_\bullet\to o$, such that $g(z_\bullet)=0$
  and each $z_d$ is a one-sided-isolated zero of $g$. (Note that the set of non-zeros of $g$ is open in $(\R^1,o).$)

   We can assume: the sequence $z_\bullet$ is  strictly
   decreasing,  and $g>0$ on some intervals $(z_d,z_d+\ep_d]$. Here $z_d+\ep_d<z_{d-1}$ for all $d$.

For each such interval construct $f_d(x)$ as in Step 3, so that:
\bee[\!\!\!\bf i.\!]
\item $f_d\in C^\infty(z_d,z_d+\ep_d)$ and $f_d$ is flat at $z_d$ and $z_d+\ep_d$;
\item  For each $d$ the expression $e^{g(x)\frac{d}{dx}}f_d(x)$ diverges on  a sequence of points $\tz^{(d)}_d\to z^+_d$.
\item For the pre-defined sequence $\{N_\bullet\}$ of equation \eqref{Eq.fast.growing.derivatives} the bound
  holds: $ \xi^{j+1}  f_d\big|_{z_d} > N_j\cdot \xi^j  f_d\big|_{z_d}$ for $j\ge d$.
  \item $\sum^{\td}_{j=0}|f^{(j)}|<\frac{1}{\td}$ on $(z_d,z_d+\ep_d)$ for all $\td\ll d.$
\eee

\

Finally, extend this collection $f_\bullet$  to a function $f\in
C^\infty(\R^1)$,
 by zeros outside of the (disjoint) intervals $\{(z_d,z_d+\ep_d)\}$.

We get: $f\in C^\infty(\R^1)$ and   $e^{g(x)\frac{d}{dx}}f(x)$
diverges on  a sequence of points $\tz^{(d)}_d\to 0^+$.
 Thus $e^{g(x)\frac{d}{dx}}f(x)\not\in C^\infty(\R^1,o)$. Moreover, for each fixed $d$ the sequence
   $\xi^{\bullet} \tf|_{\tz^{(d)}_d}$ is fast growing, as in
  \eqref{Eq.fast.growing.derivatives}.
   Thus also for the initial   derivation, $\xi\in \Der_\R(C^\infty(\R^n,o))$, we have the divergence, $e^{\xi}f \not\in C^\infty(\R^n,o)$.

\medskip

\item  Fix a derivation $\xi=\sum g_i(\ux)\frac{\di}{\di x_i}\in \Der_\R(\quots{C^\infty(\R^n,o)}{J}).$
 We construct $f\in C^\infty(\R^n,o)$ for which $e^\xi f$ diverges on a sequence of points $o\leftarrow z_\bullet\sset V(J).$
\bee[\hspace{-0.6cm}\bf\text{Part} A.]
\item  As $o\!\in\! \R^n$ is not an isolated point of $V(J),$
 we can choose a converging sequence $o\!\leftarrow \! z_\bullet\sset  V(J)$. We can assume
  that the coordinates of this sequence, $z_d\!=\!(z_{d,1},\dots,z_{d,n})$, satisfy:
  \bei
  \item $0< z_{d+1,i}<\frac{1}{2}z_{d,i}$ for any $i=1,\dots,n$ and for any $d$;
  \item $g_1(x)$ does not vanish at each $z_d$.
  \eei
Indeed, after a preliminary rotation of $\R^n,$ none of the points $z_\bullet$ lies in coordinate hyperplanes.
 Then (if needed) pass to a subsequence and permute the coordinates.

Now construct a $C^\infty$ arc $C\sset (\R^n,o)$ that contains this sequence and is locally parallel to the $x_1$-axis at each $z_j$. More precisely, we
 construct a map $[0,1)\stackrel{\ux(t)}{\to}\R^n$ satisfying:
\bei
\item $\ux(t)\in C^0[0,1)\cap C^\infty(0,1)$ and $\ux(t)=o$ \iff $t=0$.
\item $\frac{d x_1}{dt}\gneq 0$ and $\frac{d x_2}{dt}\geq 0$, \dots, $\frac{d x_n}{dt}\geq 0$ for any  $t\in (0,1)$.

Therefore the sequence $\{t_d:=\ux^{-1}(z_d)\}$ is strictly decreasing.
\item The functions $x_2(t),\dots, x_n(t)$ are locally constant at $\{t_\bullet\}.$ I.e. there exist small disjoint segments
 $\amalg(t_d-\ep_d,t_d+\ep_d)$ on which the  functions $x_2(t),\dots, x_n(t)$ are locally constant.
\eei
To construct such an arc one can start from the piecewise linear curve $\cup [z_{d+1},z_d]$ and smoothen the corners at the points $\{z_d\}$.

\item
 Define the ``restriction" of $\xi\in \Der_\R(C^\infty(\R^n,o))$ to this arc $C$ by $\de_C(t):=\frac{g_1(\ux(t))}{\frac{dx_1}{dt}}\frac{d}{dt}$.
 As $\frac{dx_1}{dt}>0$ we have $\de_C\in \Der_\R(C^\infty(0,1))$. Construct a function $f_C(t)\in C^\infty(0,1)$ such that $exp(\de_C)f_C$  diverges
  at each point $t_j$. (The construction is the same as in Steps 2,3.) Moreover, we assume the fast growing derivatives at each $z_d$,
    as in equation \eqref{Eq.fast.growing.derivatives}, and that $f$ is flat at $t=0,$ with fast vanishing derivatives.

Now ``extend" $f_C$     to $f\in C^\infty(\R^n,o)$ as follows. For each $z_d=\ux(t_d)$ and the corresponding $\ep_d$ (as in Part A) take
 the corresponding small tubular neighborhood, $\cU_d$, of the curve $\ux(t_d-\ep_d,t_d+\ep_d)$. Identify $\cU_d$ with the (trivial) bundle
  $\cU_d\stackrel{\phi_d}{\to} \ux(t_d-\ep_d,t_d+\ep_d)$. Here the fibre $\phi^{-1}_d(t)$ is a ball lying in a hyperplane orthogonal to $C$ at $\ux(t)$.

Take the $C^\infty$-bump functions, $0\le \tau_d\le 1,$ satisfying:
\beq
\tau_d(\ux)=  0  \text{ for }\ux\not\in \cU_d, \text{ and }
\tau_d(\ux)=  1,\ \text{ locally near }z_d.
\eeq
Finally define
\beq
f(\ux):=\Bigg\{\ber f_C(\phi_d(\ux))\cdot \tau_d(\ux),\ \text{ if }\ \ux\in \cU_d
\\0,\ \text{ otherwise}.
\eer
\eeq
By our assumption $f\in C^\infty(\R^n\smin o).$ Moreover, $f\in C^\infty(\R^n,o),$ as $f_C$ is flat at $o,$ with  fast-vanishing derivatives.

We have $g_1\frac{\di f(\ux)}{\di x_1}=g_1\frac{dt}{dx_1}\frac{d f_C(t)}{d t}$. And  $e^\xi f$ diverges on the sequence  $z_\bullet\sset V(J),$
 by the construction of $f_C$ (and the condition \eqref{Eq.fast.growing.derivatives}).

\eee
 Altogether, $e^\xi $ does not act on $\quots{C^\infty(\R^n,o)}{J}$.
\epr\eee

\beR
\bee[\bf i.]
\item As one sees from the proof, $e^\xi $ does not act on any subring $R\sseteq C^\infty(\R^n,o)$
 that contains enough flat functions with fast growing derivatives.
\item
   Let $R:=\quots{ C^\infty(\R^n,o)}{J}$ for an ideal  $J\supseteq \cm^\infty.$
    (This condition is close to $V(J)=0\sset (\R^n,o)$.) Then for  any $(x)$-adically nilpotent derivation $\xi\in \Der^{nilp}_\R(R)$
  the operator $e^\xi $ acts on $R$. Indeed, in this case Borel's lemma ensures the isomorphism:
\beq
\quots{C^\infty(\R^n,o)}{J}\cong \quot{C^\infty(\R^n,o)/\cm^\infty}{J\cdot \quots{C^\infty(\R^n,o)}{\cm^\infty} }
\isom{} \quot{\R[\![\ux]\!]}{J\cdot \quots{C^\infty(\R^n,o)}{\cm^\infty}}.
\eeq
  And we observe that the ring $\quot{\R[\![\ux]\!]}{J\cdot \quots{C^\infty(\R^n,o)}{\cm^\infty}}$ is $\cm$-complete. Thus any $(x)$-adically nilpotent derivation exponentiates.
\item The assumption $\xi|_{V(J)}\neq o$  is important, as the previous remark shows.
\eee
\eeR

  \end{document}